\documentclass[10pt,letterpaper]{article}

\usepackage[top=0.85in,left=2.75in,footskip=0.75in]{geometry}

\usepackage{amsmath,amssymb}

\usepackage{changepage}

\usepackage{textcomp,marvosym}

\usepackage{cite}

\usepackage{nameref,hyperref}

\usepackage[right]{lineno}

\usepackage{microtype}
\DisableLigatures[f]{encoding = *, family = * }

\usepackage[table]{xcolor}
\usepackage{multirow}
\usepackage{array}

\newcolumntype{+}{!{\vrule width 2pt}}

\newlength\savedwidth

\raggedright
\usepackage[aboveskip=1pt,labelfont=bf,labelsep=period,justification=raggedright,singlelinecheck=off]{caption}

\makeatletter
\renewcommand{\@biblabel}[1]{\quad#1.}
\makeatother

\usepackage{lastpage,fancyhdr,graphicx}
\usepackage{epstopdf}
\fancyheadoffset[L]{2.25in}
\fancyfootoffset[L]{2.25in}
\input{paper_macros.sty}
\usepackage{float}
\usepackage{subcaption}

\begin{document}
\vspace*{0.2in}

\begin{flushleft} {\Large \textbf\newline{Physiological accuracy in simulating refractory cardiac tissue: the volume-averaged bidomain model vs.~the cell-based EMI model} %
  }
\newline
\\
Joyce Reimer\textsuperscript{1},
Sebasti\'an A. Dom\'inguez-Rivera\textsuperscript{2},
Joakim Sundnes\textsuperscript{3},
Raymond J. Spiteri\textsuperscript{4*}%
\\

\bigskip \textbf{1} Division of Biomedical Engineering, University of
Saskatchewan, Saskatoon, Canada
\\
\textbf{2} Department of Mathematics and Statistics,
University of Saskatchewan, Saskatoon, Canada
\\
\textbf{3} Simula Research Laboratory, Oslo, Norway
\\
\textbf{4} Department of Computer Science, University of Saskatchewan,
Saskatoon, Canada
\\
\bigskip

* spiteri@cs.usask.ca

\end{flushleft}

\section*{Abstract}
The refractory period of cardiac tissue can be quantitatively described using strength-interval (SI) curves. The information captured in SI curves is pertinent to the design of anti-arrhythmic devices including pacemakers and implantable cardioverter defibrillators. As computational cardiac modelling becomes more prevalent, it is feasible to consider the generation of computationally derived SI curves as a supplement or precursor to curves that are experimentally derived. It is beneficial, therefore, to examine the profiles of the SI curves produced by different cardiac tissue models to determine whether some models capture the refractory period more accurately than others. In this study, we compare the unipolar SI curves of two tissue models: the current state-of-the-art bidomain model and the recently developed extracellular-membrane-intracellular (EMI) model. The EMI model's resolution of individual cell structure makes it a more detailed model than the bidomain model, which forgoes the structure of individual cardiac cells in favour of treating them homogeneously as a continuum. We find that the resulting SI curves elucidate differences between the models, including that the behaviour of the EMI model is noticeably closer to the refractory behaviour of experimental data compared to that of the bidomain model. These results hold implications for future computational pacemaker simulations and shed light on the predicted refractory properties of cardiac tissue from each model.

\section*{Author summary}
Mathematical modelling and computational simulation of cardiac activity have the potential to greatly enhance our understanding of heart function and improve the precision of cardiac medicine. The current state-of-the-art model is the bidomain model, which considers a volume average of cardiac activity. Although the bidomain model has had success in several applications, in other situations, its approach may obscure critical details of heart function. The extracellular-membrane-intracellular (EMI) model is a recently developed model of cardiac tissue that addresses this limitation. It models cardiac cells individually; therefore, it offers significantly greater physiological accuracy than bidomain simulations. This increase in accuracy comes at a higher computational cost, however. To explore the benefits of one model over the other, here we compare the performance of the bidomain and EMI models in a pacing study of cardiac tissue often employed in pacemaker design. We find that the behaviour of the EMI model is noticeably closer to experimental data than the behaviour of the bidomain model. These results hold implications for future pacemaker design and improve our understanding of the two models in relation to one another.


\section*{Introduction}
Computational cardiac models are used in both research and clinical settings for studying cardiac diseases and deriving new treatments. A commonly used framework is the bidomain model because it is thought to adequately capture the behaviour of the heart while maintaining a reasonable computational cost \cite{ref:Sundnes2006}.  The bidomain model is also often simplified to the monodomain model to further reduce the computational cost \cite{ref:Skouibine2000}. Recently, the extracellular-membrane-intracellular (EMI) model was adapted from its original form as a neuron model \cite{ref:Agudelo-Toro2013} to describe the electrical activity of the heart at the cellular level \cite{ref:Tveito2017}. Although this model poses a significantly higher computational cost than the bidomain model, it allows for the exploration of applications that require the resolution of individual cells to accurately observe cardiac function. Ongoing developments in high-performance computing make it feasible to perform meaningful simulations at such a level of detail. Due to the relative newness of the EMI model, however, there are only a handful of studies that explore its capabilities \cite{ref:Dominguez2021, ref:Tveito2017, ref:Tveito2020, Jaeger2019, Jaeger2021a, Jaeger2021d, Jaeger2021e, Jaeger2022}, and only one of these studies examines it in relation to the bidomain model \cite{Jaeger2022}. In the present study, we compare the bidomain model to the EMI model in the context of unipolar stimulation of refractory cardiac tissue, an application relevant to the design of pacemakers and defibrillators.

\begin{figure}[h]
	\centering
	\includegraphics[width=11cm]{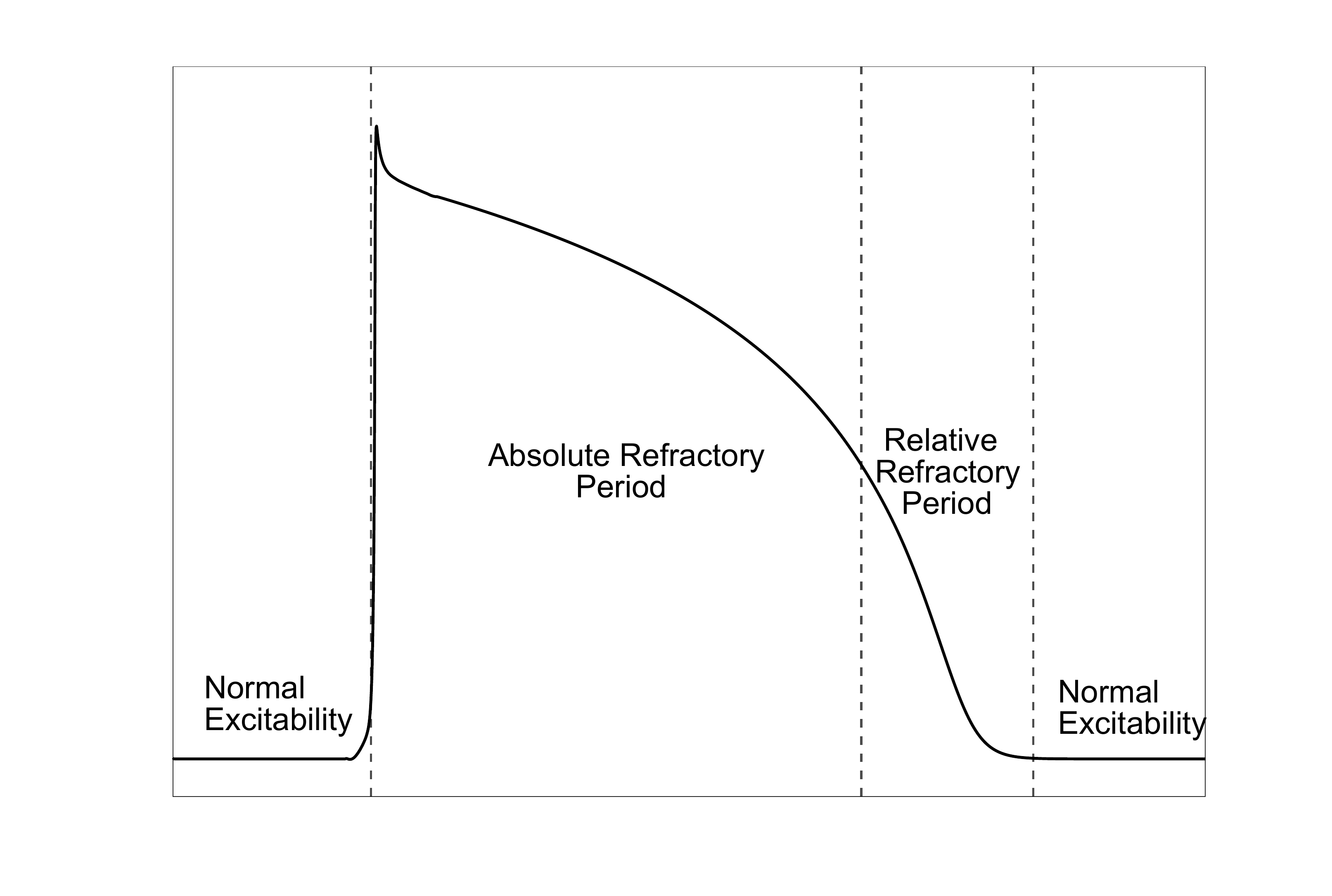}
	\caption{\textbf{The cardiac action potential and its refractory periods.} The action potential begins with an absolute refractory period, during which an action potential cannot be generated. A relative refractory period follows, in which an action potential can only be generated with a stronger-than-normal stimulus. Outside of these periods is a state of normal excitability.}\label{fig:RPs}
  \end{figure}

The refractory nature of cardiac tissue is a protective mechanism that allows for well-timed heart contractions. It is caused by the strategic closing of ion channel gates in cells that have recently fired, combined with the resulting change in tissue excitability that keeps the signal propelling in the right direction. The refractory period is subdivided into the ``absolute" and ``relative" refractory periods. In the absolute refractory period, it is not possible for another action potential to be triggered, regardless of stimulus strength. Then, after a short time, the cell enters the relative refractory period, when it becomes possible for another action potential to be initiated; however, the electrical stimulus required is higher than it was when the cell was in a resting state \cite{Boron2012, Jaye2010}. A visual representation of these stages is found in Fig.~\ref{fig:RPs}.
  
The relative refractory period is often characterized by strength-interval (SI) curves describing the strength of stimulus required to trigger another action potential at various moments during the refractory period. SI curves contain pertinent information for the design of cardiac pacemakers and implantable cardioverter-defibrillators (ICDs). These devices must deliver electrical pulses that are timed correctly and of the appropriate strength \cite{ref:Galappaththige2017a}. The two devices address opposite but related issues; therefore, a degree of precision must exist during the design process to ensure that each device is programmed to correctly address the intended issue. For example, a pacemaker is designed to correct a slow heart rate (bradycardia) \cite{ref:Reade2007}. Consequently, it needs to deliver current during the normal excitability period of the cells’ action potentials. This entails a low-energy pulse delivered shortly after the previous action potential has died out \cite{Atlee2001}. Essentially, the current delivered by a pacemaker must be high enough so as to initiate the next wave of excitation but low enough and timed correctly such that it does not trigger an arrhythmia. On the other hand, an ICD is intended to correct a fast heart rate (tachycardia) \cite{Davis2004}. Because of this, ICDs often must deliver currents while the tissue is still partially refractory \cite{Atlee2001}; i.e., while some of the cells are in the relative refractory period of the action potential. As a result, the current of an ICD must be sufficiently strong so as to override the refractoriness and successfully initiate a resynchronizing action potential. If it is not strong enough, it will be ineffective, but if it is too strong, there is a risk of damage to the cells \cite{Babbs1980}.

Each of the cardiac pacing devices described above may employ either a unipolar or bipolar electrode configuration \cite{Atlee2001, ref:Galappaththige2017a}. In the case of a unipolar configuration, a single electrode delivers either cathodal or anodal stimuli to the heart, and in the case of a bipolar configuration, two electrodes, one cathode and one anode, are spaced closely together on the heart and deliver oppositely charged electrical stimuli simultaneously. To assist with the specific design requirements of pacemakers and ICDs, SI curves are used to inform the strength of the pulse needed for the given device. 
These curves are usually generated experimentally with physical cardiac tissue (e.g., see \cite{Alferness1994, Sidorov2005}), but lately, computational versions are being explored. For example, in \cite{ref:Galappaththige2017a} and \cite{ref:Galappaththige2017b}, Galappaththige \textit{et al.} used the bidomain model to simulate unipolar and bipolar stimulation of computational cardiac tissue and plotted the resulting SI curves. 

In this study, we use the bidomain and EMI models to explore the complex interplay of tissue conductivity and membrane ion channel gating properties implicated in the cardiac refractory period. Using experimental data as reference, we calibrate the bidomain and EMI models so that they produce equivalent conduction velocities. We then use the calibrated models to find the relative refractory threshold of each model's tissue at varying time intervals after an action potential. We use these data to plot SI curves for each model and compare them against one another and against experimentally derived unipolar SI curves to determine which model gives a representation of refractory tissue that better reflects experimental data. In the resulting data, we wish to see whether defining a membrane domain in the detailed EMI model makes the EMI refractory profile more accurate than that of the spatially homogenized bidomain model.

\section*{Results}
\begin{figure}[H]
\hspace{-10mm}\includegraphics[width=15cm]{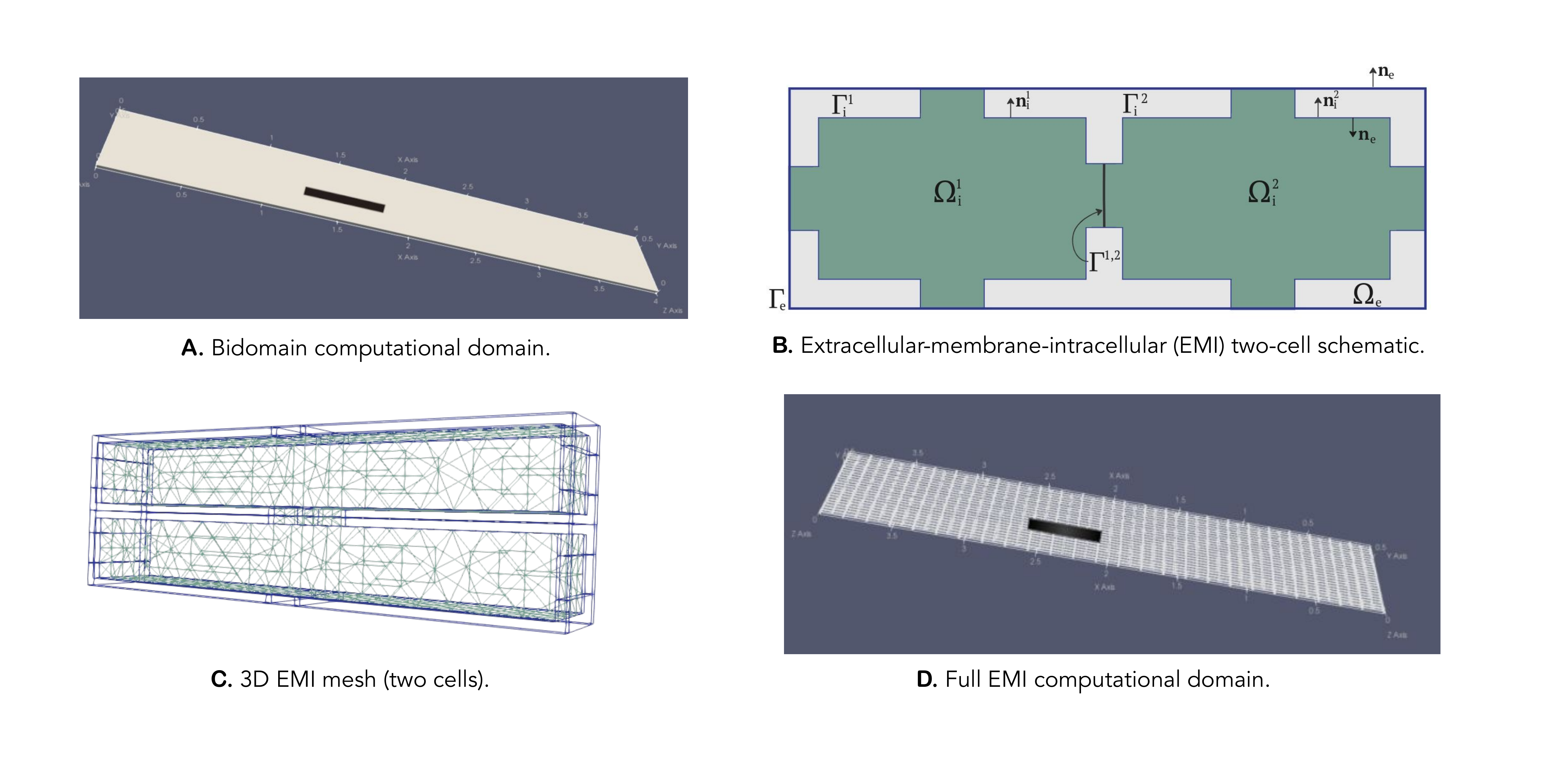}
	\caption{{\textbf{Computational domains and subdomains of the bidomain and EMI models.} \textbf{A, D}: Full computational domains of both models. The electrode position is shown in black. \textbf{B}: A schematic of two adjacent cardiac cells as defined by the EMI model. Intracellular spaces, $\Omega_\textrm{i}^1$ and $\Omega_\textrm{i}^2$, are bounded by their respective membranes, $\Gamma_\textrm{i}^1$  and $\Gamma_\textrm{i}^2$, which are in turn surrounded by extracellular space, $\Omega_\textrm{e}$. The cells are connected to one another via a gap junction, $\Gamma^\textrm{1,2}$. Outward normal vectors to the intracellular and extracellular spaces are denoted by $\mathbf{n}_\textrm{i}$ and $\mathbf{n}_\textrm{e}$. \textbf{C}: Two 3D EMI cells connected vertically via a gap junction. In green is the intracellular space, and in purple is the outline of the extracellular space.}\label{fig:compdomains}}
\end{figure}

All simulations are carried out on a three-dimensional model of rabbit ventricular cardiac tissue. The domain is a cuboid of dimension 4.0 mm $\times$ 0.625 mm $\times$ 0.025 mm (in the $x$-, $y$-, and $z$-directions, respectively). The bidomain model domain is a uniform cuboid, as shown in Fig.~\ref{fig:compdomains}A, owing to the relative simplicity of the model formulation. For the EMI simulations, however, we must define individual cardiac cells and their connectivity.  In general, each cell's intracellular space is bounded by a membrane, and all cells' membranes are surrounded by extracellular space. Adjacent cells are connected via gap junctions, which electrically couple cardiomyocytes to one another through ohmic current flow. Fig.~\ref{fig:compdomains}B shows a two-dimensional representation of this pattern with two rectangular EMI cardiomyocytes connected horizontally. An example in three dimensions is given in Fig.~\ref{fig:compdomains}C, in which two cuboid cells are coupled in a vertical configuration.  To arrive at the full computational domain shown in Fig.~\ref{fig:compdomains}D, we extend this arrangement in the $x$- and $y$-directions, resulting in a domain of $N=625$ cells in a 25 $\times$ 25 $\times$ 1 configuration. 

In the bottom third of each of the domains, we place an electrode of size 0.5 mm $\times$ 0.1 mm $\times$ 0.025 mm. A standard S1-S2 stimulus protocol is employed in all simulations. For further information on the domains and stimulus protocol, see \nameref{section:methods}.

\subsection*{The strength-interval curves}
\begin{figure}[H]
	\centering
	\includegraphics[height=9cm]{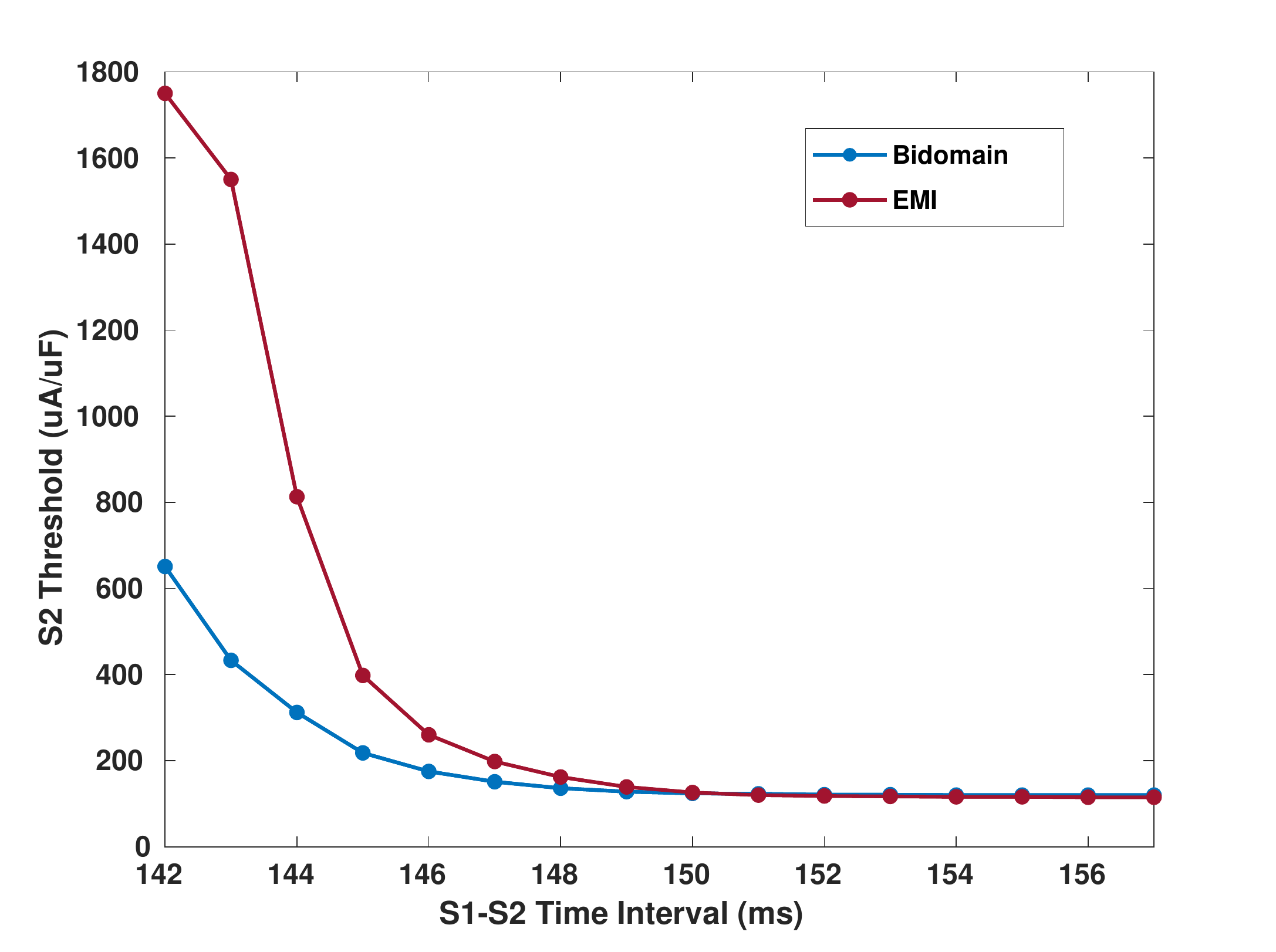}
	\caption{\textbf{The EMI model's strength-interval curve is steeper than that of the bidomain model.} These two SI curves demonstrate the different levels of refractoriness that each model predicts of cardiac tissue. The EMI model curve suggests a higher degree of refractoriness than that of the bidomain model.}\label{fig:SICurve1}
\end{figure}

Shown in Fig.~\ref{fig:SICurve1} are the bidomain and EMI SI curves acquired from our simulations in direct comparison to one another. The plots show that the conduction profiles of both models are highly similar when the S2 pulse is delivered about 149 ms or longer after the S1 pulse. However, for shorter time intervals, the curves diverge substantially, with the EMI model SI curve being much steeper than that of the bidomain model. For intervals 148 ms or shorter, the EMI simulations require a much higher S2 current than the bidomain simulations to elicit a subsequent action potential and excitatory wave. This implies that the EMI model estimates the tissue to have a higher degree of refractoriness overall relative to the bidomain model. It also suggests that the EMI model estimates the rate of change of refractoriness to be higher relative to that of the bidomain model.

\subsection*{S1-S2 experiments}
We now present two S1-S2 experiments as examples of how the SI curve in Fig.~\ref{fig:SICurve1} was generated. Here, the bidomain and EMI model S2 thresholds are determined to the nearest $\mu$A/$\mu$F for intervals from 142 ms to 157 ms. To determine whether a given S2 pulse strength produces a wave of action potentials, we take four sample points of the potential from each of the four quadrants of the domain, starting at the time of the S2 pulse. The samples are taken at the same geometric locations in each model. We analyze these traces to determine whether the S2 pulse produced a propagating wave of action potentials. It is generally quite clear from these data whether or not a wave has occurred because it shows up in an all-or-nothing manner. However, for some of the shorter time intervals, the traces reflect a more graded potential, and the decision of whether or not a propagating wave has occurred becomes more ambiguous. In these cases, we look for a peak potential of greater than 0 mV in all of the sample points. At times, there may be electrical fluctuation that brings cells close to 0 mV, but these action potentials may die out easily or have other irregular electrical properties seen in their morphology. Additional information for deciding whether or not a propagating wave has occurred is found in the time-lapse visualizations of the whole domain's changing electrical potential. We check for a situation in which the entire domain becomes depolarized by the end of the simulation, rather than only a portion.  To demonstrate the threshold determination process, action potential plots and time-lapse visualizations of the domains are shown below in Figs.~\ref{fig:aptraces}--\ref{fig:145ms-paraview} for two S1-S2 time intervals. Additional simulations are shown in Supporting Information in \nameref{S1_Fig} and \nameref{S2_Fig}.

There are two mechanisms by which excitement can occur in anisotropic cardiac tissue: ``make" and ``break". Make is when excitation of the tissue is caused directly by the depolarization of the stimulus pulse, and it begins at the onset of the pulse. Break is excitation that is caused indirectly by the relatively hyperpolarized area surrounding the electrode, and it begins once the pulse has been terminated \cite{ref:Galappaththige2017b}. Examples of each mechanism are included in the visualizations in Figs~\ref{fig:150ms-paraview} and \ref{fig:145ms-paraview}, both to illustrate the make and break mechanisms and to highlight some differences in propagation between the bidomain and EMI models.

\subsubsection*{S1-S2 interval: 150 ms}
At a time interval of 150 ms, the S2 threshold is quite apparent. To demonstrate this, the action potential plot in Fig.~\ref{fig:aptraces}A--B shows the stark difference between an S2 current strength of 123 $\mu$A/$\mu$F and 124 $\mu$A/$\mu$F for the bidomain model and between 125 $\mu$A/$\mu$F and 126 $\mu$A/$\mu$F for the EMI model.

We note here that, although the resting threshold of the EMI model is less than that of the bidomain model (115 $\mu$A/$\mu$F vs. 120 $\mu$A/$\mu$F), at this time interval, the EMI S2 threshold exceeds the bidomain S2 threshold. That is, at time intervals longer than 150 ms, the bidomain S2 threshold is higher than that of the EMI model, but at 150 ms and shorter, a crossover occurs, and the opposite is true. 

\begin{figure}[H]
	\hspace{-10mm}\includegraphics[width=15cm]{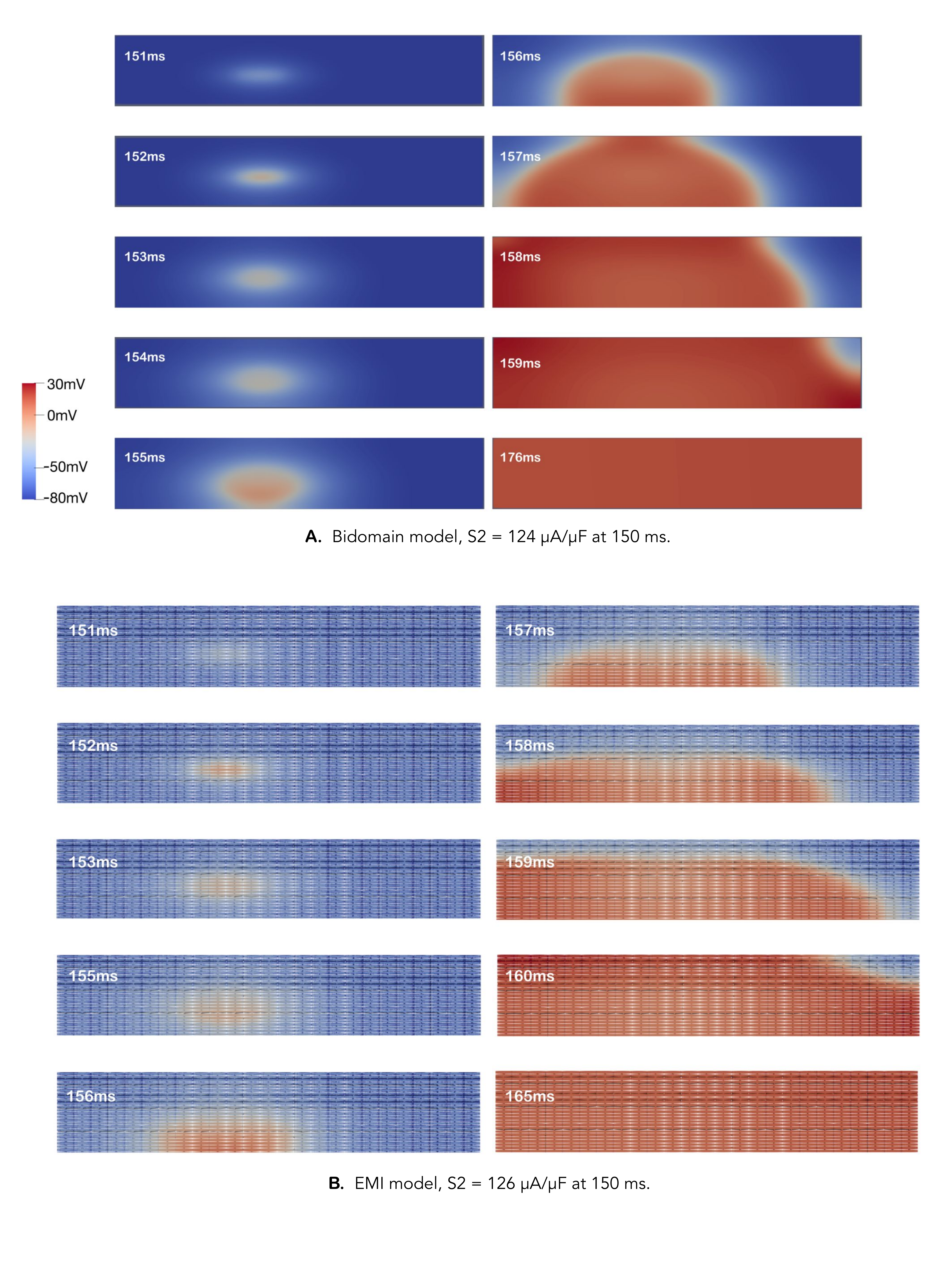} \vspace{-15mm}
\caption {\textbf{Time-lapse images illustrate propagation resulting from S2 threshold current delivered at 150 ms.} These visualizations show each model's change in transmembrane potential spreading throughout their domains over time. \textbf{A:} The bidomain S2 threshold for a 150 ms interval is 124 $\mu$A/$\mu$F. The mechanism of excitation is make, seen in the uniform outward spread of depolarization from the electrode. \textbf{B:} The EMI S2 threshold for a 150 ms interval is 126 $\mu$A/$\mu$F. The EMI model exhibits break excitation for this interval, indicated by the break of the initial depolarization centre into two propagating centres. }\label{fig:150ms-paraview}
\end{figure}

To demonstrate the conduction properties of the two models compared to one another, visualizations of the propagation produced by each model's S2 thresholds are shown in Fig.~\ref{fig:150ms-paraview}. The mechanism by which excitation occurs for each model is different at this time interval. The bidomain mechanism is make, as seen by the rapid excitation that begins at the source of the stimulus and continues to spread outward radially from that source. The EMI mechanism, on the other hand, is break, as seen in the slower conduction overall (e.g., compare frames at $t$=158 ms between \ref{fig:150ms-paraview}A and \ref{fig:150ms-paraview}B), as well as the way the excitation spreads from the part of the tissue just below the location of the electrode (in \ref{fig:150ms-paraview}B, compare frame $t$=152 ms when the electrode has just stopped firing, to $t$=156 ms when the area of excitation ``drops" and begins to propagate). Break excitation also tends to appear to split off into two circular centres that spread apart from one another, whereas make tends to appear as one circular centre that spreads outward. To this point, in \ref{fig:150ms-paraview}B, two circular centres can be identified in contrast with the single circular centre in \ref{fig:150ms-paraview}A.

\subsubsection*{S1-S2 interval: 145 ms}
\begin{figure}[H]
	\hspace{-10mm}\includegraphics[width=15cm]{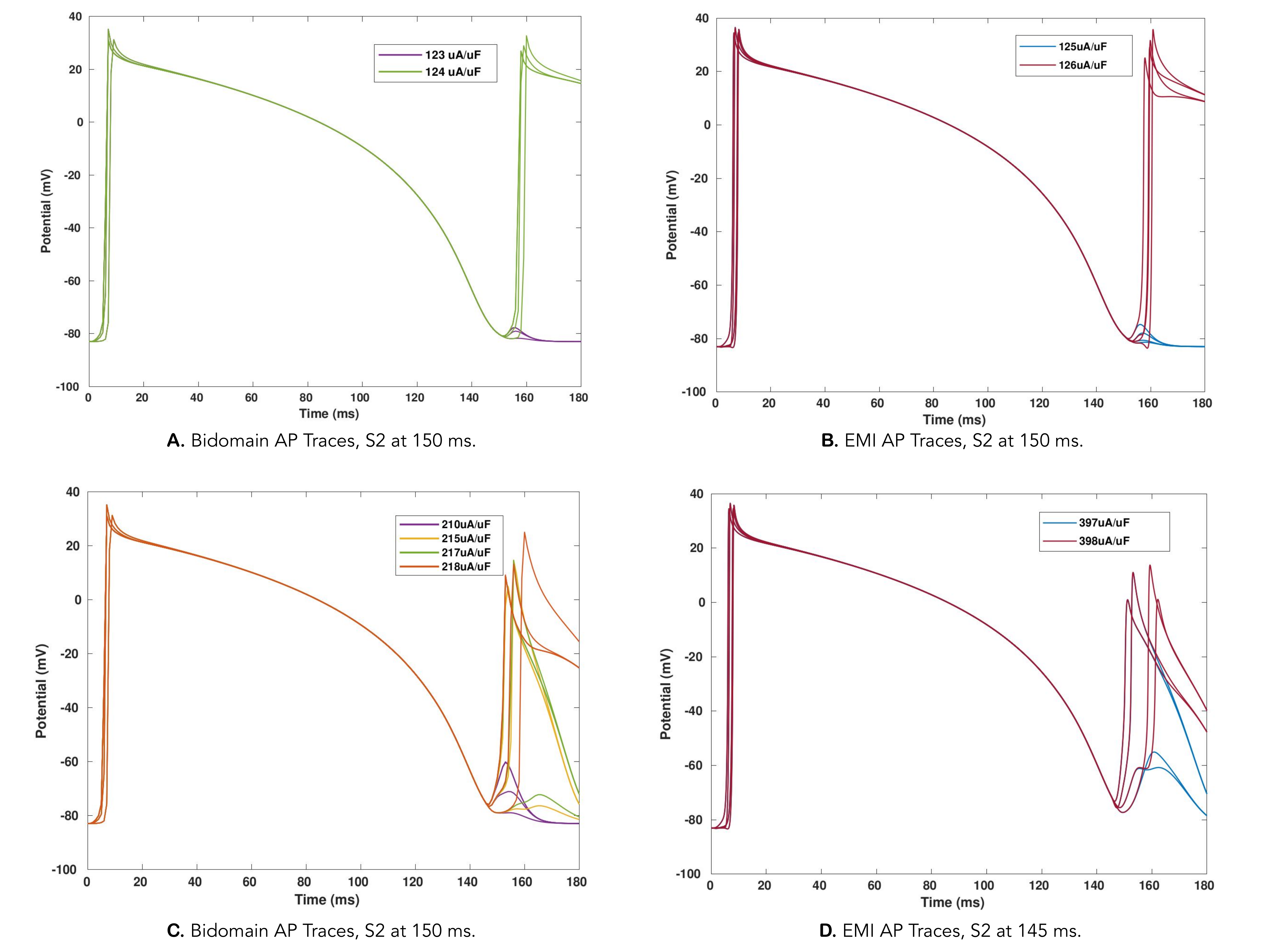}\vspace{0.75cm}
	\caption{{\textbf{Action potential traces reveal each model's S2 threshold after 150 and 145 ms.} \textbf{A}: The bidomain S2 threshold for a 150 ms interval is 124 $\mu$A/$\mu$F. \textbf{B}: The EMI S2 threshold for a 150 ms interval is 126 $\mu$A/$\mu$F. \textbf{C}: The bidomain S2 threshold for a 145 ms interval is 218 $\mu$A/$\mu$F; lower currents only trigger APs in part of the domain. \textbf{D}: The EMI S2 threshold for a 145 ms interval is 398 $\mu$A/$\mu$F.  }\label{fig:aptraces}}
\end{figure}

As an example of an interval that is less apparent, we show the plots from the simulations in which the S2 current is delivered at 145 ms. A degree of ambiguity is seen in the action potential plots as well as the visualizations. The bidomain AP traces in Fig.~\ref{fig:aptraces}C demonstrate that, for currents greater than 210 $\mu$A/$\mu$F up to 217 $\mu$A/$\mu$F, only a portion of the domain fires a second AP. At 218 $\mu$A/$\mu$F, however, the entire domain fires an AP. It is important for the S2 current to produce a wave that depolarizes the entire domain, otherwise the voltage gradient may produce a substrate for an arrhythmia. Therefore, we take 218 $\mu$A/$\mu$F to be the S2 threshold. %

The APs from the EMI simulations are similarly ambiguous. From Fig.~\ref{fig:aptraces}D, we see that 397 $\mu$A/$\mu$F produces APs only in a portion of the domain. Therefore, we choose 398 $\mu$A/$\mu$F to be the threshold because there are APs across the entire domain.

\begin{figure}[H]
	\hspace{-22mm}\includegraphics[width=16cm]{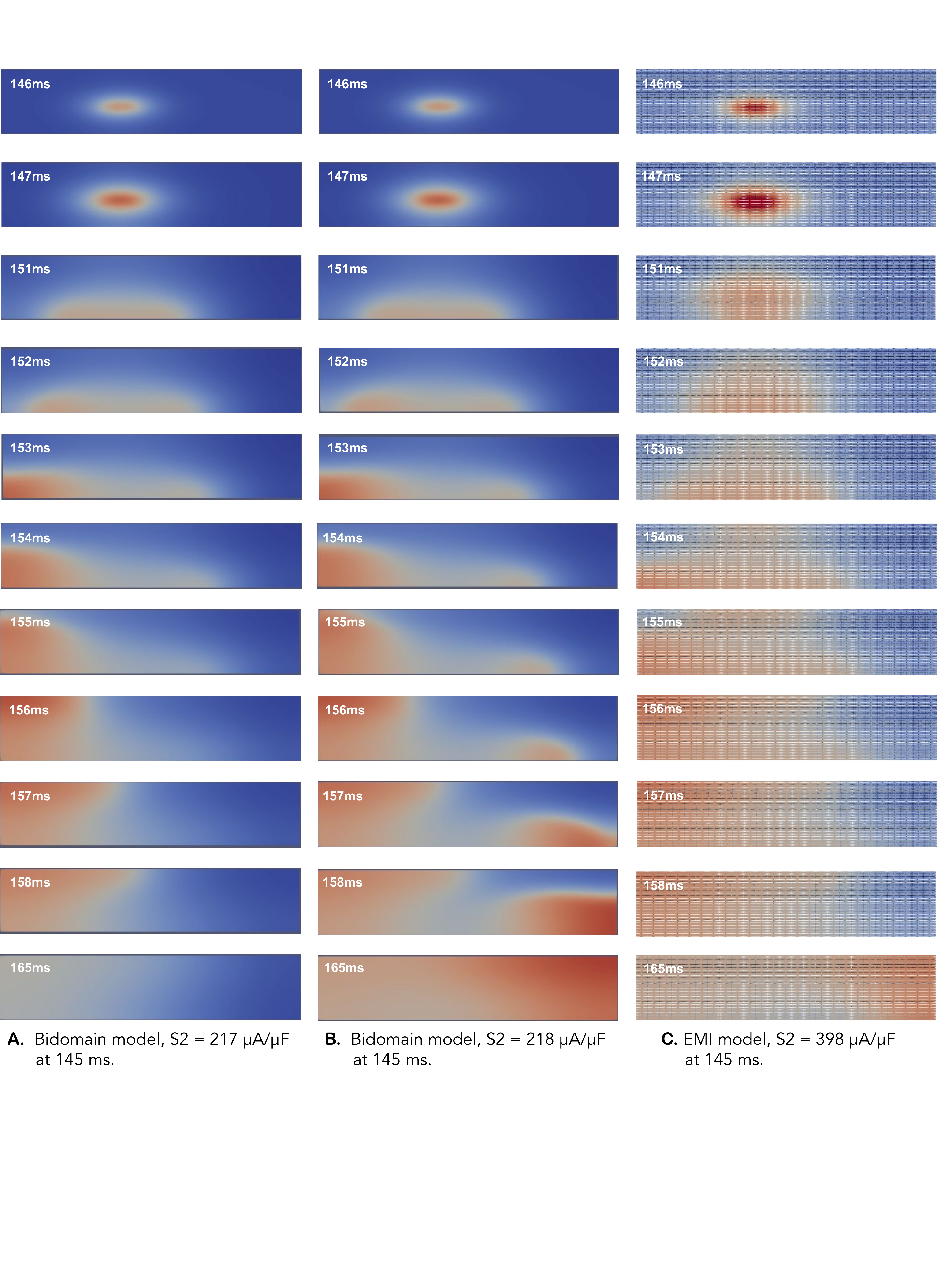} \vspace{-25mm}
\caption{\textbf{Propagation resulting from S2 currents for a 145 ms interval.} Refer to Fig.~\ref{fig:150ms-paraview} for colour scale. \textbf{A:} An S2 current of 217 $\mu$A/$\mu$F at 145 ms in the bidomain model produces an excitatory wave in only a portion of the domain that then dissipates. \textbf{B:} An S2 current of 218 $\mu$A/$\mu$F at 145 ms depolarizes the entire bidomain tissue through a break mechanism. \textbf{C:} In the EMI model, an S2 current of 398 $\mu$A/$\mu$F at 145 ms depolarizes the entire domain through a break mechanism.}\label{fig:145ms-paraview} 

\end{figure}

We can see this ambiguity reflected in the visualizations for this interval as well. For the bidomain simulations, with an S2 current of 217 $\mu$A/$\mu$F (Fig.~\ref{fig:145ms-paraview}A), there is a propagating signal that occurs; however, rather than emanating out radially via ``make" or splitting into two foci via ``break", we see the signal splitting off into two circular centres as it does in break, but then the centre on the right dissipates whereas the centre on the left continues to propagate. This results in the initiation of a spiral wave pattern, which, in physical situations, has the potential to evolve into a reentrant arrhythmia. In this simulation however, the spiral wave simply dissipates. Therefore, we choose 218 $\mu$A/$\mu$F to be the threshold (Fig.~\ref{fig:145ms-paraview}B) because it results in a more stable wave.    

We also observe this pattern in the EMI visualizations (Fig.~\ref{fig:145ms-paraview}C). Although not shown here, an S2 current of 397 $\mu$A/$\mu$F appears to set up a similarly harsh voltage gradient as was observed in the bidomain simulation at 1 $\mu$A/$\mu$F below the threshold in Fig.~\ref{fig:145ms-paraview}A. Therefore, we maintain that 398 $\mu$A/$\mu$F is the EMI threshold for 145 ms because it results in a more even depolarization of the entire domain.

\subsection*{Comparison with experimental data}
These results indicate that there is indeed a difference between the bidomain and the EMI models pertaining to the relative refractory period of the tissue. But the question remains of which refractory profile is more representative of physical cardiac tissue. To this end, we compare our computationally determined SI curves to the experimentally derived SI curves from Sidorov \textit{et al.} \cite{Sidorov2005}. Their experimental SI curves are acquired from a whole perfused rabbit heart, using unipolar cathodal stimulation with a 0.05 mm$^2$ wire electrode. There are three such curves available in the reference. Here, they are plotted and normalized by dividing the S2 threshold by the resting threshold, which for their study was 0.15 mA. The curves from this study are also normalized by their respective resting thresholds. The normalizations are done to account for the different scales of the setups. In doing this first, we are able to examine the strength-interval relationships relative to the distinct initial conditions of each setup.

\begin{figure}[h]
	\centering
	\includegraphics[width=11cm]{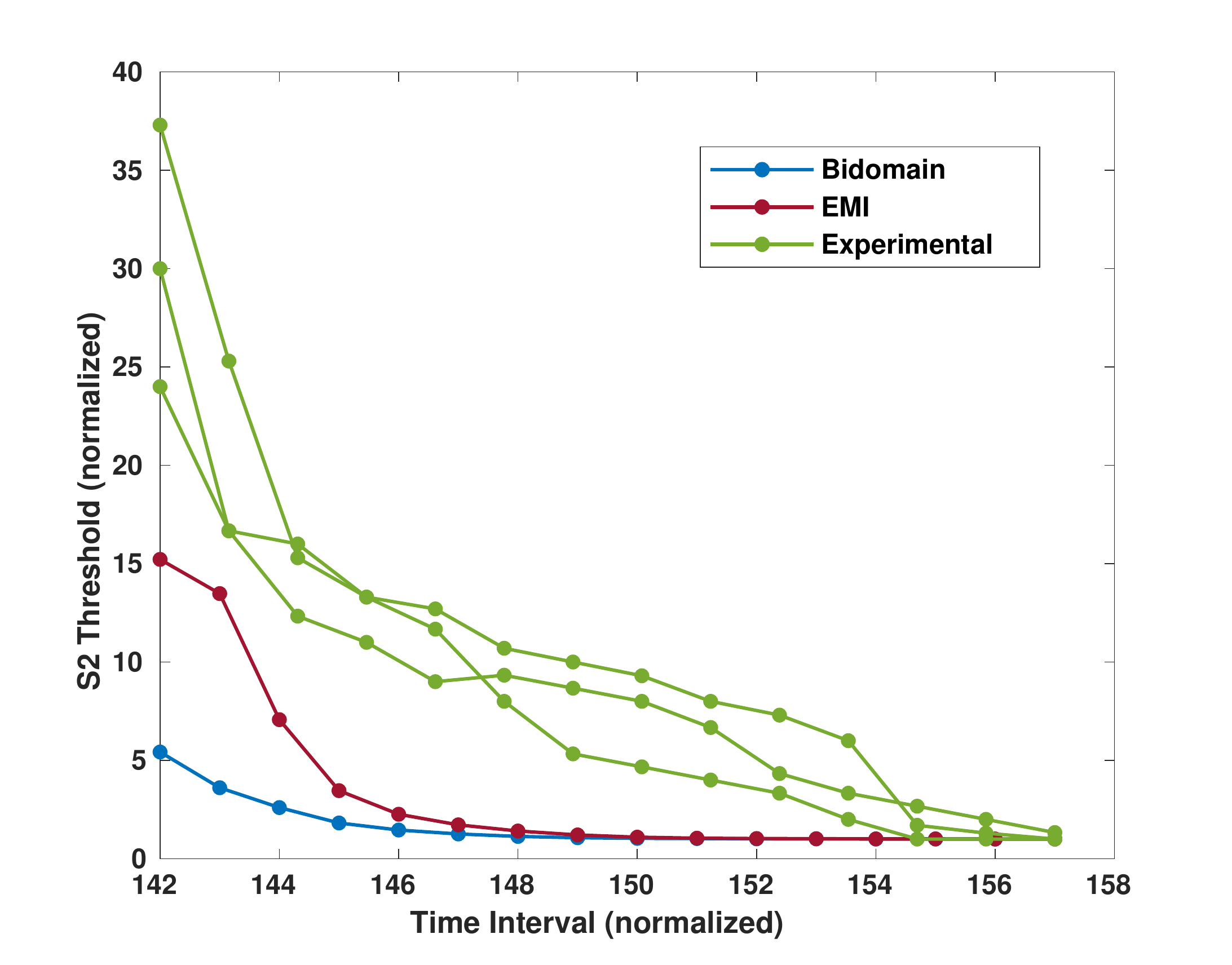}
	\caption{\textbf{The EMI model strength-interval curve is closer to experimental data than that of the bidomain model.} The EMI model produces a strength-interval curve that is closer in morphology to the experimental SI curve than that produced by the bidomain model. Each of these curves have been normalized by their respective resting thresholds.}\label{fig:SICurve2}
\end{figure}

Fig.~\ref{fig:SICurve2} shows this comparison, where we see that both models underestimate the tissue refractoriness from intervals of about 154 ms and shorter. However, we also see that the overall shape of the EMI model curve is considerably closer to the experimental curve compared to that from the bidomain model. This indicates that, given our respective implementations of each model, the EMI model depiction of the refractory period of cardiac tissue is more faithful to physiological data than that of the bidomain model.

\section*{Discussion}
We have conducted a comparative analysis between the conventional bidomain model and the contemporary EMI model. As of the time of writing, only one other comparison study between these models has been published \cite{Jaeger2022}, and the present comparison is the first to consider strength-interval curves and the refractory period of cardiac tissue. The results bring to light three notable differences between the bidomain and EMI models that can be seen in the morphology of the SI curves, the lengths of the relative refractory periods (RRPs), and the difference in their respective resting thresholds. For the cases where data were available, the EMI model displayed behaviour that was closer to experimental observations. These findings are expanded upon below.

\subsection*{Shape of strength-interval curve}\label{sec:discussion-sicurve}
As seen in Fig.~\ref{fig:SICurve1}, there is a clear difference between the shapes of the bidomain and EMI SI curves: the EMI curve is much steeper than the bidomain curve, especially at shorter S1-S2 time intervals. This indicates that the EMI model estimates the tissue to exhibit a higher degree of refractoriness during the RRP than the bidomain model. Even though tissue refractoriness is largely dependent on the gating dynamics found at the cell membrane level, which in this study is modelled in both tissues by the same cell model \cite{ref:Gray2016}, the results confirm that a component is dependent on the distribution of the cell membrane, the membrane capacitance, and other micro-discontinuities across the tissue that affect conductivity. In the bidomain model, cellular discontinuities are omitted to achieve a homogeneous domain, but in the EMI model, all cells are explicitly defined as distinct subdomains with membranes and gap junctions.

The comparison to experimental data in Fig.~\ref{fig:SICurve2} reveals that the shape of the EMI model's SI curve is closer to experimental SI curves than that of the bidomain model. This curve is normalized because the total scale of the thresholds and time intervals of the experimental data is different than the scale of our computational studies. The difference in threshold scale can be attributed to the small size of tissue modelled in comparison to the whole rabbit heart in the experimental study. It may also be attributed to the method of current delivery being different between the two setups; the experimental setup is based on a stimulus that is applied to the surface of the heart---i.e., to the extracellular domain---however, in the present study, we apply the stimulus as a transmembrane current, as if the domain is receiving a signal from neighbouring cardiac tissue. We control for these differences by normalizing each experiment's S2 thresholds by their respective resting thresholds. The reason for the discrepancy in time interval is less clear, although other bidomain studies have found a similar discrepancy (e.g., \cite{ref:Galappaththige2017a} and \cite{ref:Galappaththige2017b}). Nevertheless, these changes are controlled for by normalizing each data point by its baseline threshold and by scaling the larger time scale to the time scale of our study. These results show that, in the context of tissue refractoriness and studies involving S1-S2 stimulus protocols, the high level of detail included in the EMI model does result in an observably more accurate outcome relative to the bidomain model.

\subsection*{Length of relative refractory period}

\begin{figure}[H]
	\centering
\hspace{-.75cm}\includegraphics[width=14cm]{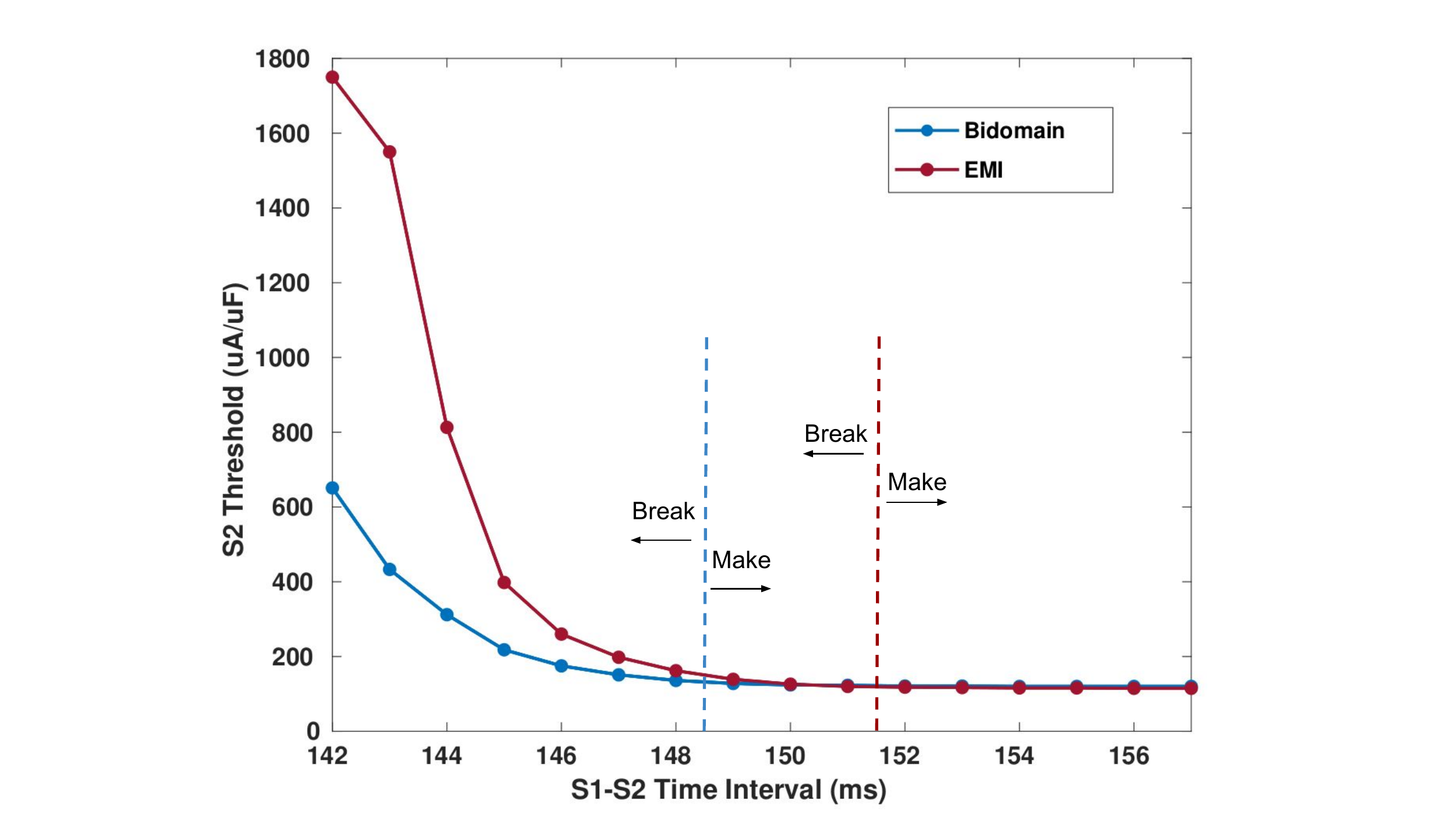}
	\caption{\textbf{The RRP of the EMI model is longer and more accurate than that of the bidomain model.} These plots show the computational SI curves from Fig.~\ref{fig:SICurve1} with make-to-break transitions marked. The bidomain model transition point is marked in blue; the EMI model transition point is marked in red. The EMI model transition happens later than the bidomain model transition, translating to a longer and more accurate RRP.}\label{fig:makebreak}
\end{figure}

In the experimental study performed by Sidorov \textit{et al.}, the end of the relative refractory period is marked by a transition from a break excitation to a make excitation that happens at around 180 ms \cite{Sidorov2005}. Then, a complete return to normal excitability happens at approximately 200 ms. In our experiments, we found that this transition point happened at a longer time interval in the EMI model compared to the bidomain model, corresponding to a longer RRP for the EMI model. Specifically, as shown in Fig.~\ref{fig:makebreak}, the transition from break to make happens in the move from 148 ms to 149 ms in the bidomain model, whereas in the EMI model, the transition is at 151 ms to 152 ms. The EMI model also had an overall longer refractory period; the bidomain model returns to a state of normal excitability at 154 ms, whereas this is 156 ms in the EMI model. Therefore, we see that the EMI model gives a more accurate representation of the length of the RRP%
.

\subsection*{Resting threshold}\label{sec:rheobase}
In this study, the EMI model has a lower resting threshold than the bidomain model. Despite calibrating the models to have similar conduction properties, the EMI model has a lower resting threshold than the bidomain model. The difference is relatively small---5 $\mu$A/$\mu$F lower. This discrepancy, however, indicates there may be structural details inherent in the EMI model that affect electrical excitability. (For details on the process of determining the resting threshold of each model, see Fig.~\ref{fig:rheobase} in \nameref{section:methods}.)

\subsection*{Limitations and future work}
The model parameters were chosen to ensure that the baseline conduction profiles of the bidomain and EMI model were identical. The cell membrane capacitances and surface-to-volume ratios of the two models were the same, and the bidomain conductivities were parameterized so that the resulting conduction velocity matched that of the EMI model. However, there may be other conductivity values that result in highly similar conduction velocities. It is also likely that other parameters besides the bidomain conductivities affect conduction velocity, such as capacitance or the EMI model conductivities. Although we limited our parameter search to the bidomain conductivities in order to preserve experimentally observed values of other parameters, it is reasonable to assume that other parameter sets exist that would give close conduction velocities in the two models. Therefore, the potential non-uniqueness of the inputs is a limitation of this study.

A second limitation is the small computational domain relative to the size of a real rabbit heart from which the experimental data were gathered. This limitation stems from the heavy computation burden required by the EMI model large domains. In light of this, an extension to this study would be an efficiently parallelized implementation of the EMI model and an examination of the effects of domain size on the SI curve. Such an examination should include the quantification of any effects the choice of domain size has on the bidomain model compared to the EMI model. Results from this type of investigation would inform tissue model choice as well as provide insight on the optimal domain size to use for each model. It should be noted that for strength-interval studies, the question of appropriate domain size is not a true parameter; rather, it is a necessary constraint to allow for solution of the EMI model. The ideal trajectory of this work is the development of software and hardware to adequately carry out whole-heart EMI simulations. Because experimental SI curves are acquired from whole hearts, only whole-heart simulations can be considered to have the ideal domain size for replicating experimental results.

Finally, the way in which the electrodes are modelled is a source of ambiguity as well. The present study uses an initial membrane current in the cell model. Physiologically, this may be thought of as endogenous excitation, as if the cells where the initial current is applied have just received the electrical signal from adjacent ventricular cells. This is how propagation occurs naturally in the heart and represents how cardiomyocytes in the centre of the heart would respond to a distally applied electrical current, whether that is an endogenous signal originating from the pacemaker cells in the sinoatrial node or from a pacemaker device providing artificial stimulation. An alternative way to model the electrode would be with the application of Neumann boundary conditions to the electrode-tissue interface. This is likely representative of the experimental setup in the physiological study of \cite{Sidorov2005}, in which an electrode is simply placed on the centre of the heart and stimulus currents are delivered through it, because an electrode placed on the exterior of the heart would deliver currents through the extracellular space rather than through adjacent intracellular space. Therefore, an area of future work would be to explore the effect of electrode modelling choice on the SI curves. Electrode modelling is an under-researched topic as a whole in cardiac simulation; although there are several ways in which electrodes are modelled in the literature (e.g., an initial transmembrane current in the cell model as was used in this study \cite{Arevalo2016, Plank2021, Prakosa2018}, Dirichlet boundary conditions representing an applied voltage \cite{Plank2021, Plank2008}, source terms in which a volume stimulus is applied \cite{Sepulveda1989, Trayanova1996}, and Neumann boundary conditions \cite{ref:Galappaththige2017a, ref:Galappaththige2017b} representing a current applied to the extracellular space), the specific methods have not been compared side by side to check for differences in their outcomes. %

\subsection*{Conclusion}
For many years, the homogenization approach of the bidomain model has been the accepted technique for simulations in the field of cardiac modelling. Work with this model has made great strides; however, its simplified approach imposes modelling limitations, which may potentially be addressed with more detailed models and greater computing power. 

In this study, we performed a comparative analysis of the bidomain and EMI models in the context of strength-interval curves and the refractory period of cardiac tissue. The aim of the analysis was to identify and explore advantages provided by the detailed but computationally expensive EMI model in relation to the less detailed but less computationally expensive bidomain model. In carrying out the study, we observed noticeable differences between the two models. These differences include the EMI model having a strength-interval curve that is closer to experimental data and a longer relative refractory period than the bidomain model. Although further analysis is needed before firm conclusions are drawn, the present findings are indicators of areas that are worth deeper exploration. Such areas include determining whether the EMI model is better suited than the bidomain model for studying refractoriness, as well as what the factors are that influence SI curves. These areas hold implications for our understanding of the refractory period of cardiac cells and for the design of pacing devices, especially in the shift toward integrating computational tools into the process of medical device development.

\section*{Methods}\label{section:methods}

\subsection*{Parameters}\label{subsection:parameters}

\subsubsection*{Calibration of conductivities}
Before carrying out the simulations, we calibrate the models so that their conduction velocities (CVs) match one another and are in agreement with available experimental data on rabbit cardiac conduction. A stimulus current is applied to the left 10\% of the domain, that is, to every point of the domain whose $x$-coordinate is between 0.0 and 0.4 mm. We then measure the longitudinal ($x$-direction) and transverse ($y$- and $z$-directions) CVs by taking five equidistant measurements across the domain along the three axes and calculating the average CV in each direction.

The conductivities in the EMI model can be viewed as absolute conductivities of the intracellular and extracellular domains, whereas the conductivities in the bidomain model are effective conductivities that are chosen to emulate both the conductivity of the tissue material and the geometry of the cells. Therefore, we first run a baseline EMI simulation with the original conductivity parameters from \cite{ref:Tveito2017} to check that the CVs acquired are within the range of experimental values. We take as experimental reference the longitudinal and transverse rabbit CVs reported in \cite{Schalij1992}, which are 61 $\pm$ 7 cm/s and 22 $\pm$ 5 cm/s, respectively.  With the baseline conductivity parameters, we found the EMI CVs to be within this range, at 61.12 cm/s and 22.08 cm/s; therefore, we use these conductivities as our final parameter values. These values, along with the other EMI parameters, are listed in Table \ref{table:emiparams}. We then optimize the bidomain conductivity tensor so that the CVs of the bidomain model match the CVs of the EMI model. 

The optimization is done using the Surrogate Optimization routine in MATLAB, which searches for the global minimum of an objective function within set bounds \cite{MatlabSurrogateOpt}. This method is designed for objective functions that are computationally expensive, such as those involving long simulations \cite{Gutmann2001}. The objective function is the squared difference between each computed bidomain CV and the corresponding experimentally validated EMI CVs, 61.12 cm/s (longitudinal) and 22.08 cm/s (transverse). We seek the bidomain conductivities that minimize this function. The starting points are set to conductivity values commonly used in the literature (e.g., \cite{ref:Sundnes2006}). The maximum number of function evaluations is set to 50. The last evaluation of the objective function showed an error of only 6.9444e$-07$ in the CVs. The resulting conductivities are listed in Table \ref{table:bdparams} and are used as the bidomain conductivity tensors. 

\subsubsection*{Tissue model parameters}

The parameters used in the bidomain and EMI simulations, including conductivities, membrane capacitance, and mesh size, are listed in Tables \ref{table:bdparams} and \ref{table:emiparams}. To ensure that the effective membrane quantity in the bidomain model is comparable to the actual membrane quantity in the EMI model, we set the value of the area of the cell membrane per unit volume ($\chi$) in the bidomain model equal to the ratio of the approximate total EMI membrane surface area to the total domain volume.  We include the width of the gap junctions in the $x$- and $y$-directions to account for the membrane contributions of the gap junctions: \begin{align*}
	\textrm{Single Cell Membrane Surface Area} &\approx (0.16 \; \textrm{mm} \times 0.025 \; \textrm{mm}) \times 2 \\
	&+ (0.16 \; \textrm{mm} \times 0.02 \textrm{mm}) \times 2 \\
	&+ (0.025 \; \textrm{mm} \times 0.02 \textrm{mm}) \times 2 \\
	&\approx 0.0154 \; \textrm{mm}^2\\
	\textrm{Total Membrane Surface Area} &\approx \; 0.0154 \textrm{mm}^2 \times 625 \; \textrm{cells}\\[4pt]
	\chi &= \frac{\displaystyle \textrm{Total Membrane Surface Area}}{\displaystyle \textrm{Domain Volume}} \; \\[4pt]
	\chi &\approx  \frac{\displaystyle 9.625 \; \textrm{mm}^2}{\displaystyle 4.0 \; \textrm{mm} \times 0.625 \; \textrm{mm} \times 0.025 \; \textrm{mm}}\\[4pt]
	\chi &\approx 154 \; \textrm{mm}^{-1}
\end{align*}
This value of $\chi$ is a slight overestimation because the gap junctions do not cover the entire width or length of the cells; therefore, we round this value down to 150 mm$^{-1}$.  We found simulations to be insensitive to $\chi$ in the neighborhood of this value.

\begin{table}[H]
	\begin{centering}
	\caption{\textbf{Bidomain parameter values.}}\label{table:bdparams}
	\begin{tabular}{|c|l|}
		\hline
		\multicolumn{1}{|l|}{\textbf{Parameter}} & \multicolumn{1}{c|}{\textbf{Value}}      \\ \hline
		$\sigma_{i,x}$                           & 0.2525 $\mu$A mV$^{-1}$mm$^{-1}$ \\ \hline
		$\sigma_{i,y}$                           & 0.0222 $\mu$A mV$^{-1}$mm$^{-1}$ \\ \hline
		$\sigma_{i,z}$                           & 0.0222 $\mu$A mV$^{-1}$mm$^{-1}$ \\ \hline
		$\sigma_{e,x}$                           & 0.821 $\mu$A mV$^{-1}$mm$^{-1}$  \\ \hline
		$\sigma_{e,y}$                            & 0.215 $\mu$A mV$^{-1}$mm$^{-1}$  \\ \hline
		$\sigma_{e,z}$                            & 0.215 $\mu$A mV$^{-1}$mm$^{-1}$  \\ \hline
		$C_m$                                    & 0.01 $\mu$F mm$^{-2}$   \\ \hline
		$\chi$                                   & 150 mm$^{-1}$        \\ \hline
		$\Delta x$                                & 0.025 mm             \\ \hline
	\end{tabular}\\
\end{centering}Conductivities ($\sigma$) are calibrated values. The proportionality constant $\chi$ is derived from the ratio of the EMI membrane surface area to the total domain volume. Membrane capacitance, $C_m$, is taken from the literature; e.g., \cite{ref:Sundnes2006}.
\end{table}

\begin{table}[H]
	\begin{centering}
	\caption{\textbf{EMI parameter values.}} \label{table:emiparams}
	\begin{tabular}{|c|l|}
		\hline
		\textbf{Parameter} & \multicolumn{1}{c|}{\textbf{Value}} \\ \hline
		$\sigma_i$         & 0.5 $\mu$A mV$^{-1}$mm$^{-1}$        \\ \hline
		$\sigma_e$         & 2.0 $\mu$A mV$^{-1}$mm$^{-1}$        \\ \hline
		$C_m$, $C^{(k,l)}$ & 0.01 $\mu$F mm$^{-2}$                \\ \hline
		$R^{(k,l)}$        & 0.15 mV$\cdot$mm${^2}\mu$A$^{-1}$       \\ \hline
		$\Delta x$         & 0.005 mm                             \\ \hline
	\end{tabular}\\
\end{centering}Adjacent cells are given labels $(k,l)$, where $ l \in \N[k]$. All parameter values are obtained from \cite{ref:Tveito2017}.
\end{table}
\subsection*{Stimulus protocol}
\begin{figure}[H] 
	\centering
	\includegraphics[height=9cm]{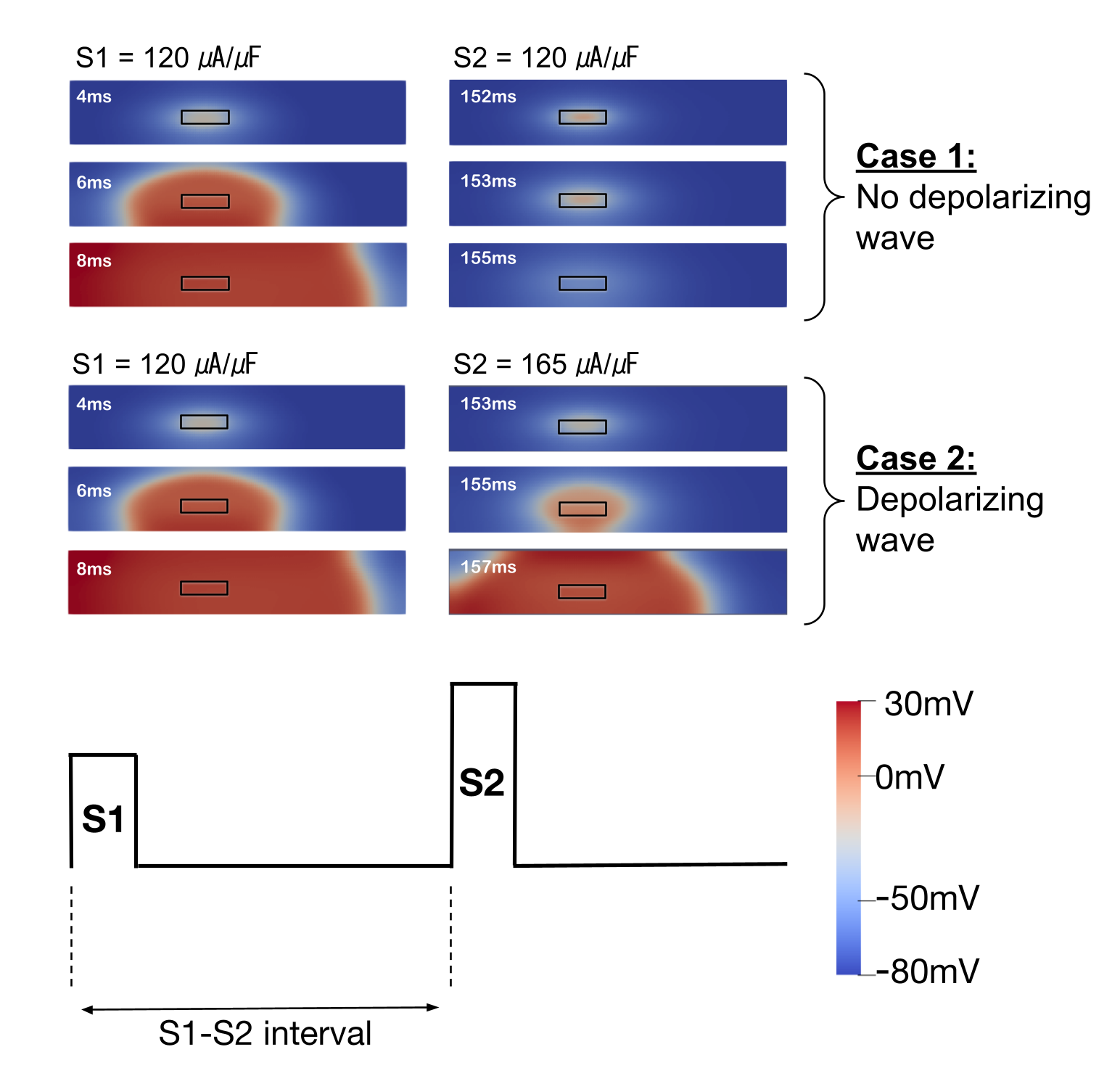}
	\caption{\textbf{The S1-S2 pacing protocol.} The position of the electrode is outlined in black. The colour scale in the bottom right indicates the transmembrane potential of the domain. The S1 current is chosen such that there is always a depolarizing wave, as seen by the spreading region of higher potential in red. The S2 threshold current is variable and thus may or may not cause a subsequent excitatory wave. In Case 1, the S2 current is too low to cause a depolarizing wave; in Case 2, the S2 current is high enough to cause a depolarizing wave. }\label{fig:protocol}
\end{figure}

The S1-S2 stimulus protocol used in this study is a common approach to pacing studies (e.g., \cite{ref:Galappaththige2017a, ref:Galappaththige2017b, Roth1996}). ``S1" refers to the first stimulus pulse that is delivered; it triggers action potentials in the cells closest to the site of stimulation that then become an excitatory wave that propagates across the entire tissue. The strength and timing of the S1 pulse is typically held constant throughout the experiment. ``S2" refers to the second stimulus pulse that is delivered after S1, and its strength and timing are variable. The aim of the protocol is to find the threshold current that is required for the S2 pulse to initiate a (depolarizing) wave of excitation for varying durations of time after the S1 pulse. At a sufficiently short time interval post-S1, the time that the S2 pulse is delivered will fall within the relative refractory period. At this point, a stronger stimulus relative to the resting threshold is required to initiate an excitatory wave. We aim to find such S2 stimulus strengths and plot the corresponding S1-S2 intervals for both the bidomain and EMI models. An illustration of this stimulus protocol is found in Fig.~\ref{fig:protocol}.

In all simulations, the S1 strength is chosen to be the threshold of current required to initiate an action potential when the cells are at rest. The S1 pulse is delivered at the beginning of the simulation ($t = 0$ ms) for a duration of 2 ms. All S2 pulses are also administered for a duration of 2 ms. The S1-S2 interval is measured as the time interval between the start times of both stimulus pulses. This implies that, because the onset of S1 is at $t = 0$ ms, the S1-S2 interval is equal to the onset time of S2. 

S2 thresholds to the nearest 1 $\mu$A/$\mu$F of intervals from 142 ms to 157 ms, incremented by 1 ms, are identified. At intervals shorter than 142 ms, the current that is required is higher than the stability of the problem allows, and at intervals longer than 157 ms, the tissue is no longer refractory in either model. The stimulus currents in this study are delivered using a cathodal unipolar electrode. The primary reason for this choice of electrode is the availability of physiological data on cathodal unipolar SI curves \cite{Sidorov2005} to which we may compare our results.
The electrode dimensions are also chosen such that the surface area is equal to the electrode surface area in contact with the tissue in \cite{Sidorov2005} with the aim of a close comparison.

\subsection*{Determination of resting thresholds} \label{section:threshold}

\begin{figure}[H]
	\hspace{-10mm}\includegraphics[width=15cm]{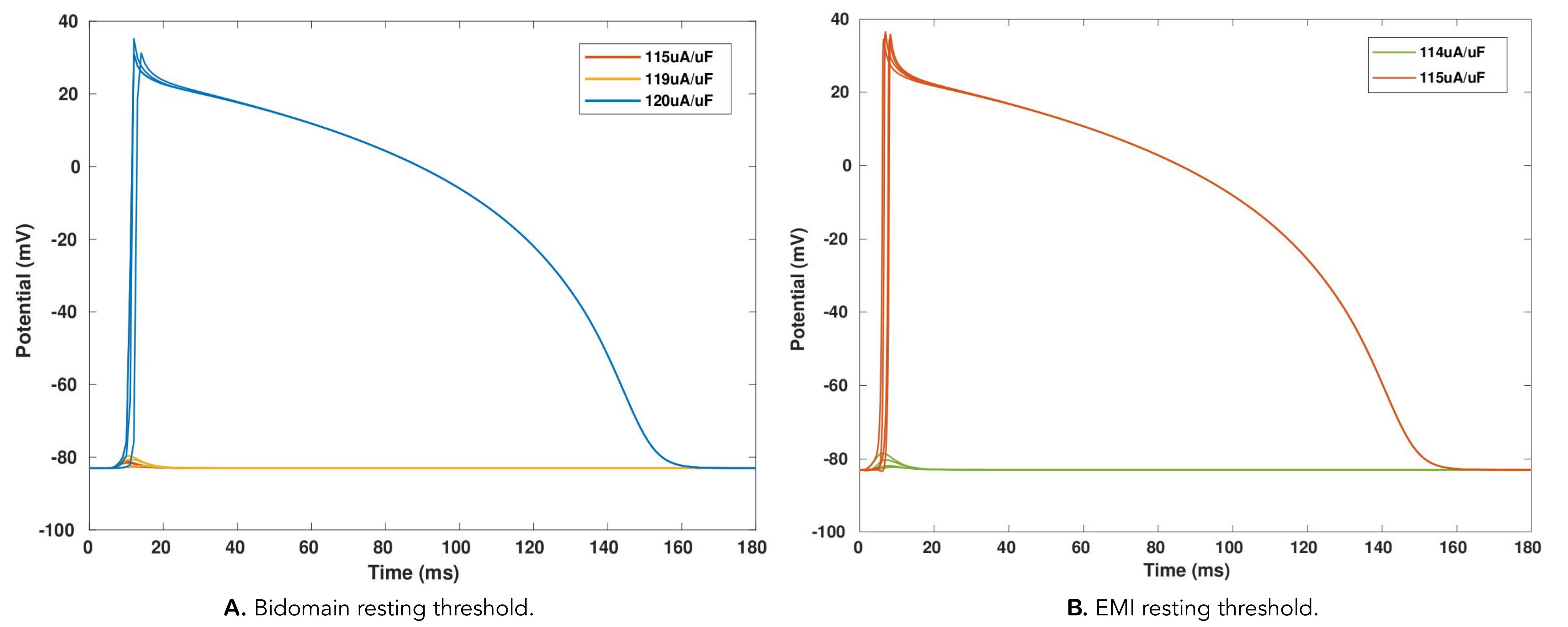}\vspace{0.75cm}
	\caption{{\textbf{Determination of the bidomain and EMI model S1 thresholds.} The S1 threshold value is the lowest stimulus current that is able to trigger an action potential from rest. \textbf{A}: The bidomain threshold is 120 $\mu$A/$\mu$F. \textbf{B}: The EMI threshold is 115 $\mu$A/$\mu$F.   }\label{fig:rheobase}}
\end{figure}

To determine whether a given S1 pulse strength produces a wave of action potentials, four sample points of the potential are taken at different locations on the membrane, plotted as a function of time, and analyzed. These sample points are taken at the same geometrical location (near the four corners of the domain) for both bidomain and EMI simulations. To within 1 $\mu$A/$\mu$F, the thresholds were found to be 120 $\mu$A/$\mu$F for the bidomain model and 115 $\mu$A/$\mu$F for the EMI model. These values are used as the S1 pulse strengths for all experiments in the study. To illustrate, the plots in Fig. \ref{fig:rheobase} show the threshold current triggering action potentials on each domains. %

\subsection*{Bidomain model}\label{subsection:bdmodel}
The bidomain model is a mathematical description of the electric activity of excitable cells in living tissue from a macroscopic point of view. Specifically, it couples a system of linear diffusion partial differential equations describing the propagation of electrical activation through the tissue, typically referred to as the tissue model, with a system of non-linear reaction partial differential equations describing the flow of ions across the membrane of individual cells in the tissue, referred to as the cell model \cite{ref:Cervi2018}. No explicit description of the geometry of the cells or any component of the structure of the tissue is given in this model, rather, the intra- and extracellular spaces and membrane components are each located everywhere in space. Despite the non-physiological homogenization of the subcellular spaces, the bidomain model is currently considered the gold standard for describing the electrical activity in the heart (e.g., \cite{Gillette2021}).

Here, we present a mathematical description of the bidomain model. Let $\Omega$ be a bounded domain in $\rrr^3$, representing three-dimensional cardiac tissue, and let $\tn[0],\tend>0$ be such that $\tn[0]<\tend$. On the spatial domain $\Omega$ and the time interval $[\tn[0],\tend]$, the bidomain model is formulated as
\begin{subequations}\label{eq:bidomain-eqs}
	\begin{align}
		\pderiv[]{\s}{t} &= \f(t,\s,v),\label{eq:bidomain-cellmodel}\\
		\chi \Cm\pderiv[]{v}{t} +\chi (\Iion(\s,v) + \Istim) &= \nabla\cdot(\sigmai\nabla v)+\nabla\cdot(\sigmai\nabla \ue),\label{eq:bidomain-eq}\\
		0 &= \nabla\cdot(\sigmai\nabla v) + \nabla\cdot((\sigmai+\sigmae)\nabla\ue)\label{eq:bidomain-constraint},
	\end{align}
\end{subequations}
subject to the following boundary conditions on $\partial\Omega\times[\tn[0],\tend]$ given by
\begin{subequations}\label{eq:bidomian-bc}
	\begin{align}
		(\sigmai\nabla v + \sigmai\nabla \ue)\cdot \hatn &= 0,\label{eq:bidomian-zerobc}\\
		(\sigmae\nabla \ue)\cdot\hatn &= 0.\label{eq:bidomian-istimbc}
	\end{align}
\end{subequations}
The vector field $\f$ represents the non-linear cell model, which is itself a function of time $t\in[\tn[0],\tend]$, the state vector field $\s=\s(\x,t)$, describing the cellular state at the location $\x\in\Omega$, and the transmembrane potential $v = v(\x,t)$. The ionic current, $\Iion$, is a non-linear term derived from the cell model, and the S1 or S2 stimulus current, $\Istim$, is applied to the cell model within the region of the electrode. The unknowns of the system include $v$ as well as the extracellular potential,  $\ue = \ue(\x,t)$. The parameters $\sigmai$ and $\sigmae$ are tensors that represent the conductivities of the intracellular and extracellular spaces, respectively. The parameter $\Cm$ is the capacitance of the cell membrane per unit area, $\chi$ is the area of the cell membrane per unit volume, and $\hatn$ is the outward unit normal vector to $\partial\Omega$. For all bidomain simulations in this study, $\Omega$ is a cuboid of dimension 4.0 mm $\times$ 0.625 mm $\times$ 0.025 mm (in the $x$-, $y$-, and $z$-directions, respectively).

\subsection*{EMI model}\label{subsection:EMImodel}
The EMI model also describes the behaviour of excitable cells in living tissue. In contrast to the bidomain model, the EMI model explicitly considers three physical domains: the intracellular space, the cell membrane, and the extracellular space. The tissue is described as a collection of cells distributed across the extracellular space that interact with each other through capacitative current and through gap junctions in their membranes. In this sense, the EMI model represents a more detailed description of the structure of living tissue at the cellular level than the bidomain model \cite{ref:Tveito2020}.

We now give a detailed mathematical description of the EMI model. Let $\Omega$ be a bounded domain in $\rrr^3$, representing three-dimensional cardiac tissue. We assume that $N$ three-dimensional cardiac cells are embedded in $\Omega$. For $k = 1,2,\ldots,N$, we denote cell $k$ in $\Omega$ by $\Omegai[k]$ and its boundary by $\partial\Omegai[k]$. The extracellular space is defined as the three-dimensional space outside these cardiac cells and inside the domain $\Omega$. We denote the extracellular space by $\Omegae$ and its boundary by $\partial\Omegae$. If $\Omegai[k]$ and $\Omegai[l]$ are two adjacent cells, i.e., $l \in \N[k]$, $k,\ l=1,2,\ldots,N$, where $\N[k]$ is the \textit{adjacency set} of $k$, we define the gap junction between these two cells by $\GammaGap[k,l] = \partial\Omegai[k]\cap\partial\Omegai[l]$. Similarly, the membrane of the cell $\Omegai[k]$, $k=1,2,\ldots,N$, is defined by $\Gammai[k]:=\partial\Omegai[k]\cap\partial\Omegae$. The outer boundary of the domain $\Omega$ is denoted by $\Gammae$.

Although an arbitrary number of cells with complex geometries may be represented with the EMI model, the EMI simulation domain in this study is comprised of $N=625$ individual cuboid cardiomyocytes arranged in a 25 $\times$ 25 $\times$ 1 configuration. The main cuboid that forms one cardiac cell is 0.155 mm long in the $x$-direction, with a height and width of 0.020 mm in the $y$- and $z$-directions, comparable to cell dimensions observed in \cite{Wiegerinck2006}. Adjacent cells are spaced 0.005 mm apart in both the $x$- and $y$-directions and are connected to one another via smaller cuboids of size 0.005 mm $\times$ 0.010 mm $\times$ 0.005 mm. These smaller cuboids represent the gap junctions through which neighbouring cells are electrically coupled. A box comprised of extracellular material of size 0.16 mm $\times$ 0.025 mm $\times$ 0.025 mm is built around each cell. %

The EMI model gives the extracellular and intracellular potentials in the cardiac tissue, as well as the potentials across each cell membrane and gap junction. The extracellular potential is denoted by $\ue$, and the intracellular potential of cell $k$ is denoted by $\ui[k]$, $k=1,2,\ldots,N$. The transmembrane potential $\vk[k]$ across the intracellular-extracellular membrane $\Gammai[k]$, $k=1,2,\ldots,N$, is defined as %
  \begin{equation*}
  \vk[k] = \ui[k] - \ue.
\end{equation*}
The potential $\wk[k,l]$ across gap junction $\GammaGap[k,l]$ is defined by
\begin{equation*}
	\wk[k,l] = \ui[k] - \ui[l] ,\ l \in \N[k],\ k=1,2,\ldots,N.
\end{equation*}
The EMI model is formulated as follows. On the spatial domains
$\Omegae$ and
$\Omegai[k]$, $k=1,2,\ldots,N$, the extracellular and intracellular potentials satisfy 

\begin{subequations}\label{eq:emi-constraints}
	\begin{align}
      0 &= \nabla\cdot(\sigma_e\nabla\ue), &\textrm{in}\ &\Omegae,
\label{eq:emi-extracell-constraint}      \\
		0 &= \nabla\cdot(\sigma_i\nabla\ui[k]), &\textrm{in}\ &\Omegai[k],\label{eq:emi-intracell-constraint} \\
	0 &= (\sigma_e\nabla \ue)\cdot\hatn, & \textrm{on}\ &\Gammae,
	\end{align} 
  while on the membrane $\Gammai[k]$, the cellular state vector $\sk[k](\x,t)$ and transmembrane potential $\vk[k]$ satisfy 
	\begin{align}
		\frac{\partial\sk[k]}{\partial t} &= \fk[k](t,\sk[k],\vk[k]),&\quad\text{on $\Gammai[k]$},\label{eq:emi-cellmodel}\\
		\Cmk[k]\pderiv[]{\vk[k]}{t} &= \Imk[k] - \Iionk[k] - I_{\textrm{stim}}^{(k)},&\text{on $\Gammai[k]$},\label{eq:emi-ap-cellmembrane}\\
		\Imk[k] = (\sigma_e\nabla\ue)\cdot\nex &= -(\sigma_i\nabla\ui[k])\cdot\nin[k],&\text{on $\Gammai[k]$},\label{eq:emi-current-cellmembrane}
	\end{align}
and finally on the gap junctions $\GammaGap[k,l]$, the transjunction potential
$\wk[k,l]$ satisfies
  \begin{align}
		\Cmk[k,l]\pderiv[]{\wk[k,l]}{t} &= \Imk[k,l] - \Igap[k,l],&\text{on $\GammaGap[k,l]$},\label{eq:emi-ap-gapjunction}\\
		\Imk[k,l] = (\sigma_i\nabla\ui[k])\cdot\nin[k] &= -(\sigma_i\nabla\ui[l])\cdot\nin[l],&\text{on $\GammaGap[k,l]$}.\label{eq:emi-current-gapjunction}    
  \end{align}
\end{subequations}

The ion current $\Iionk[k]$ and the vector field $\fk[k]$ are non-linear terms related to the cell model. The stimulus current $I_{\textrm{stim}}^{(k)}$ (S1 or S2) is known and applied on $\Gammai[k]$ to the cell model within the region of the electrode. The parameter $\Cmk[k]$ is the capacitance of the membrane $\Gammai[k]$ per unit area, the ion current $\Imk[k]$ is the steady-state current density across the membrane $\Gammai[k]$, and $\nin[k]$ is the outward unit normal vector to $\partial\Omegai[k]$. For $l \in \N[k]$, $k=1,2,\dots,N$, the parameter $\Cmk[k,l]$ is the capacitance of the gap junction $\GammaGap[k,l]$, $\Imk[k,l]$ is the steady-state ion current density across the gap junction $\GammaGap[k,l]$, and $\Igap[k,l]=\Igap[k,l](\wk[k,l])$ is the ion current due to the potential $\wk[k,l]$ across the gap junction. We note that at the gap junction $\GammaGap[k,l]$, the normal vectors $\nin[k]$ from $\partial\Omegai[k]$ and $\nin[l]$ from $\partial\Omegai[l]$ satisfy $\nin[k] = -\nin[l]$. The parameters $\sigma_i$ and $\sigma_e$ are constants representing the conductivities of the intracellular and extracellular spaces, respectively.

\subsection*{Cell model}
Both the bidomain and the EMI models are coupled to a cell model via the terms $\Iion$ and $\Iionk[k]$, respectively. The dynamics of such cell models are governed by a system of ODEs. Following other S1-S2 computational studies \cite{ref:Galappaththige2017a} and \cite{ref:Galappaththige2017b}, we use the Gray--Pathmanathan parsimonious model of rabbit myocyte action potentials \cite{ref:Gray2016}. This cell model has 13 parameters, the values of which are given in \nameref{S1_Table}.

The cell model represents active ion current flow, but the EMI model also includes a passive mechanism of ohmic current flow $\GammaGap[k,l]$ at the gap junctions $\Igap[k,l]$, defined by \begin{align}
	\Igap[k,l] = \frac{1}{\Rgap[k,l]}\wk[k,l],
\end{align}
where $\wk[k,l]$ is the potential across the gap junction and $\Rgap[k,l]$ is the resistance of the gap junction, $k=1,2,\ldots,N$, $l \in \N[k]$.

\subsection*{Temporal and spatial discretization}\label{subsection:implementation}
The non-linear PDE systems describing the bidomain and EMI models are discretized in time and space for their numerical solution. Operator splitting is used to separate the solution of the cell model from that of the tissue model to take advantage of specialized solvers for each part of the ensuing split model. A spatial discretization is performed by the finite element method, yielding a set of ODEs for the cell model and linear set of differential-algebraic equations (DAEs) for the tissue model as per the method of lines.  A basic multi-rate method is then used to advance the time integration of the ODEs and DAEs. Details of this discretization are now given.

\subsubsection*{Temporal discretization: operator splitting plus multi-rate integration}
The time evolution of the PDE system describing the bidomain model as well as the EMI model is performed using operator splitting \cite{ref:Cervi2018,ref:Tveito2020}. Operator splitting is a common divide-and-conquer approach to solving differential equations. The idea is that a given problem can be divided into sub-problems that are easier to handle than the entire problem. The result is an algorithm that is either more tractable or efficient than one applied to the problem as a whole.

We start by applying operator splitting to the bidomain model. To this end, we first gather the unknowns $v$, $\ue$, and $\s$ into a vector $\yBD = (v,\ue,\s)^\transpose$ and rewrite the equations in \eqref{eq:bidomain-eqs} and \eqref{eq:bidomian-bc} as the PDE system
\begin{align}\label{eq:bidomain-daesystem}
	\pderiv[]{\aa\yBD}{t} = \fimBD(\yBD) + \fexBD(\yBD),\quad \yBD(\tn[0]) = {\yBD}_0,
\end{align}
where ${\yBD}_0$ is a given initial data vector and $\aa$ is rectangular matrix of coefficients such that $\aa\y = (\s,v,0,0,0)^\transpose$. The functions $\fimBD(\y)$ and $\fexBD(\y)$ are defined by 
\begin{align*}
	\fimBD(\y) = \left(\begin{matrix}
		\zero\\
		\frac{1}{\chi\Cm}\big(\nabla\cdot(\sigmai\nabla v) + \nabla\cdot(\sigmai\nabla\ue)\big)\\
		\nabla\cdot(\sigmai\nabla v) + \nabla\cdot((\sigmai+\sigmae)\nabla\ue)\\
		(\sigmai\nabla v+\sigmai\nabla\ue)\cdot\hatn\\
		\sigmae\nabla\ue\cdot\hatn %
	\end{matrix}\right)
\end{align*}
and 
\begin{align*}
	\fexBD(\y) = \left(\begin{matrix}
		\f(t,\s,v)\\
		-\frac{1}{\Cm}(\Iion(\s,v) + \Istim)\\
		0\\
		0\\
		0
	\end{matrix}\right).
\end{align*}
Splitting the equations describing the bidomain model in this way isolates the linear tissue model from the non-linear cell model. This requires only linear solves from implicit time-stepping methods applied to the tissue model and allows us to apply explicit time-stepping methods to evolve the cell variables. The notation $\fim$ and $\fex$ is meant to foreshadow this choice of numerical methods to be applied to the individual operators. Although common, such a choice is by no means the only one nor is it necessarily optimal.

A similar operator splitting is applied to the EMI model. If we gather all unknowns in the EMI model in the vector $\z = (\ui[k],\ue,\s)^\transpose$, then the equations in \eqref{eq:emi-constraints} can be rewritten as 
\begin{align}\label{eq:emi-daesystem}
	\pderiv[]{\bb\z}{t} = \gim(\z) + \gex(\z),\quad\z(\tn[0]) = {\yEMI}_0,
\end{align}
where $\zn[0]$ is a given initial data vector and the rectangular matrix $\bb$ is chosen to be such that $\bb\z = (0,0,0,\s,\v,0,\w,0)^\transpose$. The functions $\gim(\y)$ and $\gex(\y)$ take the form
\begin{align*}
	\gim(\y)=\left(\begin{matrix}
		\nabla\cdot(\sigma_e\nabla\ue)\\
		\nabla\cdot(\sigma_i\nabla\ui[k])\\
		(\sigma_e\nabla\ue)\cdot\hatn \\
		\zero\\
		-\frac{1}{\Cmk[k]}(\sigma_i\nabla\ui[k])\cdot\nin[k]\\
		(\sigma_i\nabla\ui[k])\cdot\nin[k] + (\sigma_e\nabla\ue)\cdot\nex\\
		-\frac{1}{\Cmk[k,l]}\big((\sigma_i\nabla\ui[l])\cdot\nin[l]+\Imk[k,l]\big)\\
		(\sigma_i\nabla\ui[k])\cdot\nin[k] + (\sigma_i\nabla\ui[l])\cdot\nin[l]
	\end{matrix}\right),
\end{align*}
and
\begin{align*}
	\gex(\y)=\left(\begin{matrix}
		0\\
		0\\
		0\\
		\fk[k](t,\sk[k],\vk[k])\\
		-\frac{1}{\Cmk[k]}(\Iionk[k] + I_\textrm{stim}^{(k)})\\
		0\\
		0\\
		0
	\end{matrix}\right),
\end{align*}
where we have left the indices $\l \in \N[k]$, $k = 1,2,\dots,N$, in to denote the block structure of the entries
and where we have replaced the values of $\Imk[k] = -(\sigma_i\nabla\ui[k])\cdot\nin[k]$ (cf. \eqref{eq:emi-current-cellmembrane}) on $\Gammai[k]$ and $\Imk[k,l] = -(\sigma_i\nabla\ui[l])\cdot\nin[l]$ (cf. \eqref{eq:emi-current-gapjunction}) on $\GammaGap[k,l]$.

Let $\cc$ be a rectangular matrix either denoting the matrix $\aa$ or $\bb$, and let $\y$ be a vector denoting either the vector $\yBD$ or $\z$. Then the systems in \eqref{eq:bidomain-daesystem} and \eqref{eq:emi-daesystem} can be re-written as
\begin{align}\label{eq:daesystem}
	\pderiv[]{\cc\y}{t} = \him(\y) + \hex(\y),\quad\y(\tn[0]) = {\y}_0,
\end{align}
where the functions $\him$ and $\hex$ denote either the functions $\fimBD$ and $\fexBD$ or $\gim$ and $\gex$, respectively. Given a fixed time step size $\dt$, we apply the first-order Godunov operator splitting to advance the state of the system in \eqref{eq:daesystem} from $\tn[n]$ to $\tn[n+1] = t + \dt$:
\begin{enumerate}
	\item Solve $\pderiv[]{\cc\tilde\y}{t} = \hex(\tilde\y)$ with $\tilde\y(\tn[n]) = \y(\tn[n])$ to find $\tilde\y(\tn[n+1])$.
	\item Solve $\pderiv[]{\cc\hat\y}{t} = \him(\hat\y)$ with $\hat\y(\tn[n]) = \tilde\y(\tn[n+1])$ to find $\hat\y(\tn[n+1])$. 
	\item Assign $\y(\tn[n+1]) = \hat\y(\tn[n+1])$. 
\end{enumerate}
All simulations in this study use $\Delta t = 0.1$ ms.  The Godunov method introduces a first-order splitting error; we also use first-order time integrators for Steps 1 and 2 in a multi-rate fashion to ensure stability. Specifically, the Forward Euler (FE) method with $M=100$ equal time sub-steps is used for Step 1 to find $\tqn[n+1] = \tqn[n+1,M]$ as an approximation for $\tilde\y(\tn[n+1])$. 
 \begin{align}\label{eq:femethod}
	\cc\tqn[n,j+1] = \cc\tqn[n,j] + \frac{\dt}{M}\hex(\tqn[n,j]),\quad j=0,1,\ldots,M-1,\quad \tqn[n,0] = \y(\tn[n]),
\end{align}
In situations when the applied stimulus current, $\Istim$, is large, a more stable method than FE is required to find $\tilde\q(\tn[n+1])$. In these cases, we use the Rush--Larsen (RL) method, a strategy based on partitioning the cell model ODEs into stiff and non-stiff components \cite{Marsh2012}. In particular, the gating variables tend to be stiff, so these are split off and solved with an exponential integrator, whereas the non-gating variables that are less stiff are solved with FE. We apply this method to all gating variables within the given cell model and continue to use FE in the form of Eq. \ref{eq:femethod} for the non-gating variables. To perform Step 2, we use the Backward Euler method,
\begin{align}\label{eq:bemethod}
	\cc\hqn[n+1] - \dt\,\him(\hqn[n+1]) = \cc\tqn[n+1],
\end{align}
to compute $\hqn[n+1] \approx\y(\tn[n+1])$. 

\subsubsection*{Spatial discretization: finite element method}
The finite element method is used to construct the linear system \eqref{eq:bemethod}. using P1c-P1c finite elements for approximating the transmembrane potential, $v$, and the extracellular potential, $\ue$, in the bidomain model, where P1c denotes the Lagrange elements of degree 1 on the triangle with inter-element continuity. In order to discretize the spatial variables of the bidomain model in \eqref{eq:bemethod}, we first integrate by parts the terms involved in this equation to obtain a weak formulation (in space) of the system. Let $w,\we\in H^1(\Omega)$ be two test functions. Then, integrating by parts \eqref{eq:bemethod} and using the Neumann boundary condition for the intracellular and extracellular potentials (cf. \eqref{eq:bidomian-bc}), we obtain the system
\begin{align}\label{eq:weakbemethodbidomain}
    \frac{1}{\dt}\aa\hyn[n+1] + \fimBD(\hyn[n+1]) = \frac{1}{\dt}\aa\tyn[n+1],%
\end{align}
where the matrix $\aa$ is such that
\begin{align*}
    \aa\hyn[n+1] = \left(\begin{matrix}
        \displaystyle\int_\Omega\vn[n+1]w\ dV\\
        0
    \end{matrix}\right),\quad 
    \aa\tyn[n+1]=\left(\begin{matrix}
        \displaystyle\int_\Omega\tvn[n+1]w \ dV\\
        0
    \end{matrix}\right),
\end{align*}
the vector $\fimBD(\y)$ is defined as
\begin{align*}
    \fimBD(\y) = \left(\begin{matrix}
        \displaystyle\int_\Omega\frac{\sigmai}{\chi\Cm}\Big(\nabla v + \nabla\ue\Big)\cdot\nabla w\ dV\\
        \displaystyle\int_\Omega\frac{1}{\chi\Cm}\Big(\sigmai\nabla v + (\sigmai+\sigmae)\nabla\ue\Big)\cdot\nabla\we \ dV
    \end{matrix}\right).
\end{align*}

We now describe the details of the spatial discretization of the EMI model. Following a similar approach to that in \cite{ref:Dominguez2021,ref:Tveito2017}, to construct the linear system for the EMI model, where $\cc = \bb$, $\q = \z$, and $\h = \g$, we employ a mixed RT0-P0-P1c scheme to find approximations to the ion current in $\Omega$, the intracellular and extracellular potentials $\ui[k]$ and $\ue$, the transmembrane potential $\vk[k]$, and the potential across the gap junction, $\wk[k,l]$. Here, RT0 denotes the Raviart–Thomas finite element of order zero, P0 denotes piecewise constant discontinuous finite elements, and P1c are the continuous Lagrange elements of degree 1, as defined above. 

For completeness, we give details of the discrete formulation. We define $\Gammam = \bigcup_{k=1}^N\Gammai[k]$ and $\Gammag=\bigcup_{\substack{k,l=1\\k\neq l}}^N\GammaGap[k,l]$ to denote the cellular membranes and gap junctions in $\Omega$. To discretize the space variables, we first define the flux $\q\in {H(\nabla\cdot;\Omega)}$, the potential $u\in L^{2}(\Omega)$, the transmembrane potential $\vm\in \honehalf[\Gammam]$, and the gap transmembrane potential $\wg\in\honehalf[\Gammag]$ by \begin{align}
  \q&=\begin{cases}
    \sigma_i\nabla \ui[k],&\quad\text{in $\Omegai[k]$},\quad k = 1,2,\ldots,N,\\
    \sigma_e\nabla \ue,&\quad\text{in $\Omegae$},
    \end{cases}\label{eq:qdef}\\
    u&=\begin{cases}
    \ui[k],& \quad\quad\quad\,\,\text{in $\Omegai[k]$},\quad k = 1,2,\ldots,N,\\
    \ue,& \quad\quad\quad\,\,\text{in $\Omegae$},
    \end{cases}\label{eq:udef}
\end{align}
and
\begin{align} \label{eq:vdef}
\vm|_{\Gammai[k]} = \vk[k],\quad \wg|_{\GammaGap[k,l]}=\wk[k,l],\quad k= 1,2,\ldots,N,\quad l \in \N[k].
\end{align}
Let $\p\in H(\nabla\cdot;\Omega)$ be a test function, and let {$\hatn$ be the normal vector on $\Gammam\cup\Gammag$, chosen to} point outward from the intracellular domains $\Omegai[k]$, $k=1,2,\ldots,N$, into the extracellular domain $\Omegae$, and from $\Omegai[k]$ to
$\Omegai[l]$, $k\neq l$. Also, denote by $\sigma$, $\Cm$, $\Cg$ as non-zero parameters defined by
$\sigma|_{\Omegai[k]} = \sigma_i$, $\sigma|_{\Omegae} = \sigma_e$, $\Cm|_{\Gammai[k]} = \Cmk[k]$, and $\Cg|_{\GammaGap[k,l]}=\Cmk[k,l]$, for $k=1,2,\ldots,N$, $l \in \N[k]$. Then, integrating
by parts \eqref{eq:bemethod}, we obtain 
\begin{align*}
\int_\Omega\frac{1}{\sigma}\q\cdot\p \ dV =& \sum_{k=1}^N\int_{\Omegai[k]}\nabla \ui[k]\cdot\p \ dV+ \int_{\Omegae}\nabla \ue\cdot\p \ dV\\
=& -\int_{\Omega}u\nabla\cdot\,\p \ dV +\int_{\Gammam}\vm\,\p\cdot\hatn \ dS +\int_{\Gammag}\wg\,\p\cdot\hatn \ dS.
\end{align*}
Moreover, for $r\in L^{2}(\Omega)$, we obtain
{\begin{align*}
 \int_\Omega r\nabla\cdot\,\q \ dV =& \sum_{k=1}^N\int_{\Omegai[k]} r\nabla\cdot(\sigma_i\nabla\ui[k]) \ dV + \int_{\Omegae} r\nabla\cdot(\sigma_e\nabla\ue) \ dV = 0.
\end{align*}}
Finally, for $s\in\honehalf[\Gammam]$, we have
\begin{align*}
	-\int_{\Gammam}s\,\q\cdot\hatn \ dS + \frac{\Cm}{\dt}\int_{\Gammam}\vmn[n+1]s \ dS = \frac{\Cm}{\dt}\int_{\Gammam}\tvmn[n+1]s \ dS,
\end{align*}
and for $z\in\honehalf[\Gammag]$, we have 
\begin{align*}
	- \int_{\Gammag}z\,\q\cdot\hatn \ dS+ \frac{\Cg}{\dt}\int_{\Gammag}\Big(1+\frac{1}{\Rg}\Big)\wgn[n+1]z  \ dS = \frac{\Cg}{\dt}\int_{\Gammag}\wgn[n]z \ dS.
\end{align*}
To follow a similar structure to that in \eqref{eq:weakbemethodbidomain}, we can collect the previous four equations to form the following system
\begin{align}\label{eq:weakbemethodemi}
	\frac{1}{\dt}\dd\hzn[n+1] + \gim(\hzn[n+1]) = \frac{1}{\dt}\dd\tzn[n+1],%
\end{align}
where the matrix $\dd$ is such that
\begin{align*}
    \dd\hzn[n+1] = \left(\begin{matrix}
				0\\
				0\\
        \Cm\displaystyle\int_{\Gammam}\vmn[n+1]s \ dS\\
				\Cg\displaystyle\int_{\Gammag}\Big(1+\frac{1}{\Rg}\Big)\wgn[n+1]z \ dS
    \end{matrix}\right),\quad 
    \dd\tzn[n+1]=\left(\begin{matrix}
			0\\
			0\\
			\Cm\displaystyle\int_{\Gammam}\tvmn[n+1]s \ dS \\ \Cg\displaystyle\int_{\Gammag}\wgn[n]z \ dS
    \end{matrix}\right),
\end{align*}
and the vector $\gim$ is defined as
\begin{align*}
    \gim(\z) = \left(\begin{matrix}
			\displaystyle\int_\Omega\frac{1}{\sigma}\q\cdot\p \ dV + \int_{\Omega}u\nabla\cdot\,\p \ dV - \int_{\Gammam}\vm\,\p\cdot\hatn \ dS - \int_{\Gammag}\wg\,\p\cdot\hatn \ dS\\
			\displaystyle\int_\Omega r\nabla\cdot\,\q \ dV\\
			\displaystyle-\int_{\Gammam}s\,\q\cdot\hatn \ dS\\
			\displaystyle-\int_{\Gammag}z\,\q\cdot\hatn  \ dS
    \end{matrix}\right).
\end{align*}

The systems in \eqref{eq:weakbemethodbidomain} and \eqref{eq:weakbemethodemi} are weak formulations of the systems in \eqref{eq:bemethod} for the bidomain and the EMI models, respectively. The state variables $\s$ and $\sk[k]$, $k=1,\ldots,N$, do not appear in \eqref{eq:weakbemethodbidomain} and \eqref{eq:weakbemethodemi} because of how the operator splitting is defined in Step 1 and Step 2 of the previous section. Otherwise, given that $\f$ and $\fk[k]$ are non-linear in $\s$ (and $\sk[k]$), a non-linear system of equations would need to be solved instead.

\subsection*{Software implementation}

The bidomain experiments are run using a modified version of the Cancer, Heart and Soft Tissue Environment (Chaste)~\cite{Chaste}. Chaste is a high-performance simulation software package that can be used for various biological modelling applications, including cardiac electrophysiology, tumour cell division, cell population behaviour, and lung ventilation. In terms of cardiac modelling, Chaste has the capacity to perform bidomain and monodomain simulations. It is considered a high-quality code because of its use of modern programming techniques and robust self-verification measures, including various test suites embedded in the structure of the code \cite{ChastePLOS,Chaste}. It has been independently verified in multiple instances across the literature (e.g., \cite{ChastePath,Niederer2011}).  

The default solution method for the bidomain model in Chaste is the semi-implicit method \cite{Chaste}, which solves the bidomain equations as a coupled system rather than splitting the system into explicit and implicit components. However, because of the ability for the Godunov method to solve the bidomain model with greater efficiency than the semi-implicit method \cite{Spiteri2016}, we modify the latest version of Chaste to allow us to employ Godunov splitting in the experiments. This choice also allows us to use the same solution method in the EMI model without introducing an additional confound. 

An open-source version of the EMI model is not available; therefore, we built an implementation using FreeFem++~\cite{ref:Hecht2012}. FreeFem++ is a software library written in C++ that is designed for the numerical solution of PDEs using the finite element method. To solve the linear systems posed by the EMI model, the FreeFem++ code interfaces with the library MUMPS (Multifrontal Massively Parallel Sparse Direct Solver) \cite{MUMPS:2,MUMPS:1}. MUMPS is a direct solver for sparse matrices%
. In this study, the linear system~\eqref{eq:weakbemethodemi}
is sparse and symmetric positive definite, allowing for the use of the sparse Cholesky factorization.

The cuboid domain used in this study has a relatively simple geometry overall. Because of its simplicity, we are able to use the meshing tools built into Chaste for the bidomain model. However, due to the high level of detail implicit in the EMI model, we use the mesh generator software Gmsh \cite{Gmsh} outside of FreeFem++ to develop a script that can generate a mesh with an arbitrary number of rectangular cells in the $x$-, $y$-, and $z$-directions. 

\section*{Supporting information}

\paragraph*{S1 Table.}
\label{S1_Table} 
{\bf Cell model parameters.} Parameter values for parsimonious rabbit action potential model. Values taken from \cite{ref:Galappaththige2017a, ref:Gray2016}.
\begin{table}[H]
	{
		\centering
		\begin{tabular}{|l|l|l|}
			\hline
			\multicolumn{1}{|c|}{\textbf{Parameter}} & \multicolumn{1}{c|}{\textbf{Definition}}                                                               & \multicolumn{1}{c|}{\textbf{Value}} \\ \hline
			$V_{\textrm{rest}}$ & Resting membrane potential      & $-83$ mV          \\ \hline
			$g_{Na}^{\max}$     & Maximal conductance, $I_{Na}$   & 0.11 mS/mm${^2}$ \\ \hline
			$E_{Na}$           & Reversal potential, sodium      & 65 mV           \\ \hline
			$E_{K}$            & Reversal potential, potassium   & $-83$ mV          \\ \hline
			$E_{h}$                         & \begin{tabular}[c]{@{}l@{}}Sodium channel half-inactivated\\ potential\end{tabular}           & $-74.7$ mV                    \\ \hline
			$E_{m}$                         & \begin{tabular}[c]{@{}l@{}}Sodium channel half-activated\\ potential\end{tabular}             & $-41$ mV                      \\ \hline
			$k_m$              & Sodium channel activation slope & 4 mV            \\ \hline
			$k_r$              & Parameter for shape of AP       & 21.28 mV        \\ \hline
			$k_h$                           & \begin{tabular}[c]{@{}l@{}}Sodium channel inactivation\\ slope\end{tabular}                   & 4.4 mV                      \\ \hline
			$\tau_m$           & Sodium activation time constant & 0.12 ms         \\ \hline
			$\tau_{ho}$                     & \begin{tabular}[c]{@{}l@{}}Sodium inactivation time \\ constant scale\end{tabular}            & 6.80738 ms                  \\ \hline
			$\delta_h$                      & \begin{tabular}[c]{@{}l@{}}Sodium channel inactivation time\\ constant asymmetry\end{tabular} & 0.799163                   \\ \hline
			$g_K$              & Conductance of $I_K$ at $E_K$   & 0.003 mS/mm$^2$ \\ \hline
		\end{tabular}\par
	}\end{table}
\vspace{6cm}
\paragraph*{S1 Fig.}
\label{S1_Fig} 
{\bf Failed excitation at 150 ms.} These visualizations show what failed excitation looks like in each model, at 1 $\mu$A/$\mu$F below the S2 threshold for an S1-S2 interval of 150 ms. (Refer to Fig.~\ref{fig:150ms-paraview} for colour scale.)
\includegraphics[height=19cm]{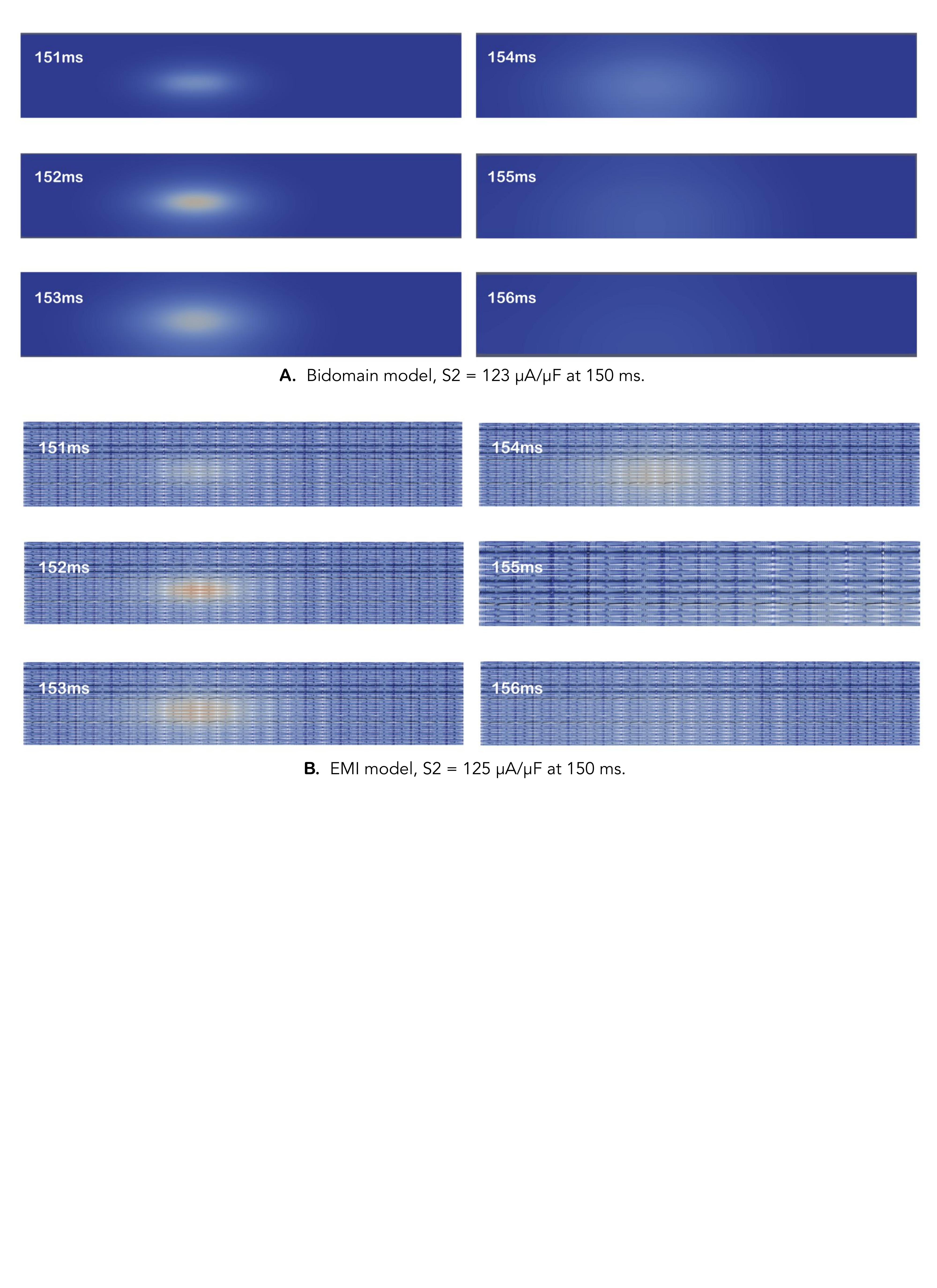}\\
	\label{fig:150ms-fail-paraview}
\vspace{3cm}
\paragraph*{S2 Fig.}
\label{S2_Fig}
{\bf AP plots and visualizations for S1-S2 Interval 148 ms.} As an additional example of the determination of the S2 threshold, the action potential plots and visualizations for the 148 ms S1-S2 interval are shown here. \textbf{A:} The bidomain S2 threshold for a 148 ms interval is 136 $\mu$A/$\mu$F. \textbf{B:} The EMI S2 threshold for a 148 ms interval is 162 $\mu$A/$\mu$F. \textbf{C:} Propagation pattern in bidomain model, S2 = 136 $\mu$A/$\mu$F. Break excitation. \textbf{D:} Propagation pattern in EMI model, S2 = 162 $\mu$A/$\mu$F at 148 ms. Break excitation.
	\begin{center}\hspace{-35mm} \includegraphics[height=11cm]{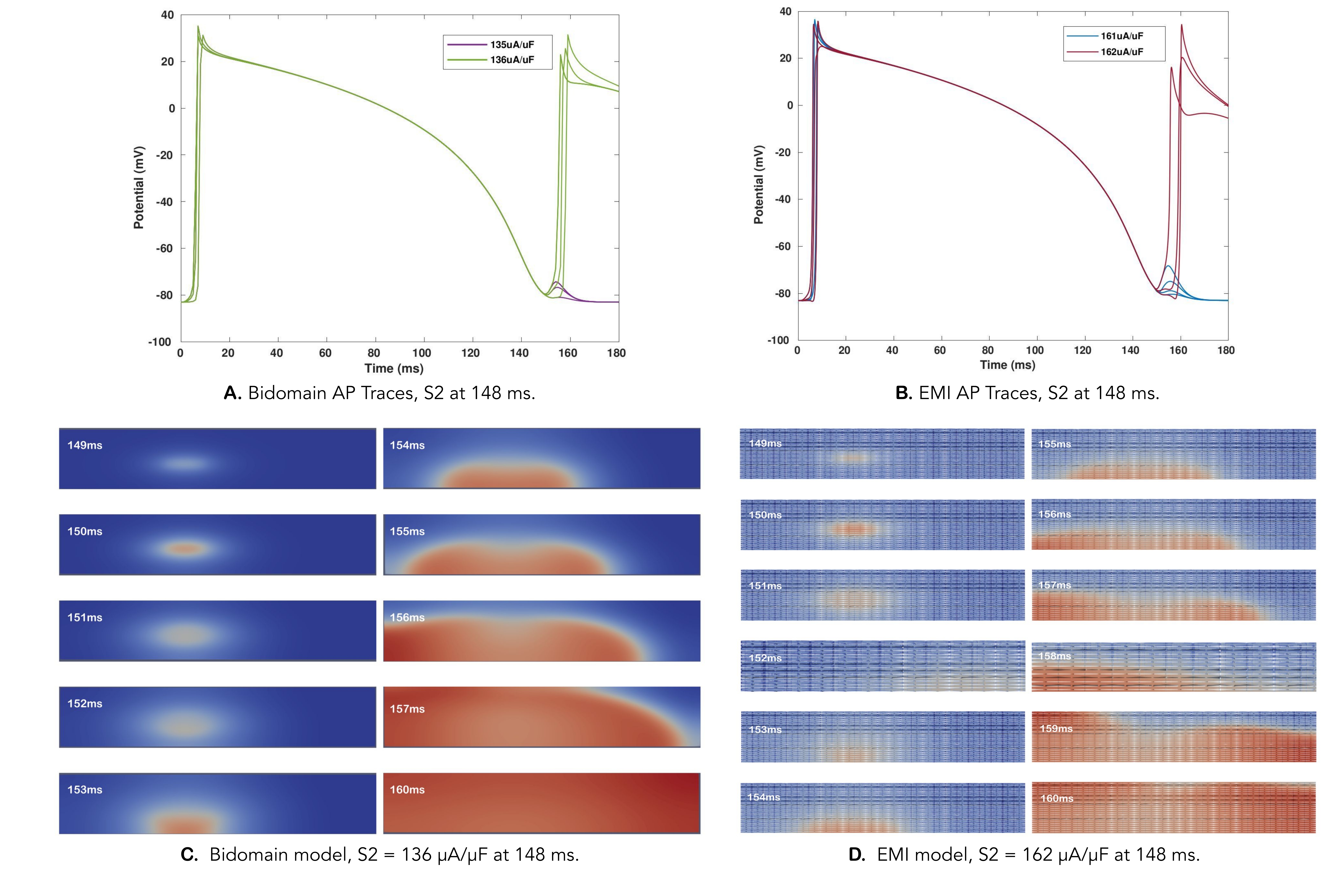}\\\end{center}

\section*{Acknowledgments}

The authors acknowledge support from the Natural Sciences and
Engineering Research Council of Canada through its Discovery Grant
Program (RGPN-2020-04467) and its Postgraduate Scholarships and the
Pacific Institute for the Mathematical Sciences through its
postdoctoral fellowship program.

\nolinenumbers

%
%
%
\bibliography{references}{}

\end{document}